\documentclass{article}
\usepackage[a4paper, left=2cm, right=2cm, top=2.5cm, bottom=1.5cm]{geometry}
\usepackage{subfiles}
\usepackage{mathtools}
\usepackage{amssymb}
\usepackage{graphicx,color}
\usepackage{amsthm}
\usepackage{amsmath}
\usepackage{enumerate}
\usepackage{tikz}
\usepackage{subcaption}
\usepackage{float}
\usepackage{subfig}
\usepackage{rotating}
\usepackage{array}
\usepackage{xfrac}

\newtheorem{theorem}{Theorem}
\newtheorem{prop}{Proposition}
\newtheorem{corollary}{Corollary}

\newtheorem{remark}{Remark}
\newtheorem{lemma}{Lemma}

\newtheorem{example}{Example}
\newtheorem{assumption}{Assumption}

\def\ba{\begin{array}}
	\def\ea{\end{array}}
\def\bi{\begin{itemize}}
	\def\ei{\end{itemize}}
\def\mR{\mathbb{R}}
\def\mZ{\mathbb{Z}}

\def\mE{\mathbb{E}}
\def\mP{\mathbb{P}}
\def\m1{1}

\def\cG{\mathcal{G}}
\def\cC{\mathcal{C}}
\def\eps{\varepsilon}

\providecommand{\Sp}{\mathrm{Sp}}
\newcommand{\AG}[1]{{A_{#1}^\cG}}

\let\oldhat\hat

\renewcommand{\hat}[1]{\oldhat{\mathbf{#1}}}
\def\tr{\mathrm{tr}\hspace{1mm}}

\usepackage[normalem]{ulem}

\begin{document}
\title{On convergence rate bounds for Markovian push-sum algorithms\thanks{The research was supported by NRDI (National Research, Development and Innovation Office) grant KKP 137490.}
}

\author{Bal\'azs Gerencs\'er\thanks{HUN-REN Alfr\'ed R\'enyi Institute of Mathematics, Budapest, Hungary and E\"otv\"os Lor\'and University, Department of Probability and Statistics, Budapest, Hungary, (email: gerencser.balazs@renyi.hu)} \and%
Mikl\'os Kornyik\thanks{HUN-REN Alfr\'ed R\'enyi Institute of Mathematics, Budapest, Hungary, (email: kornyik.miklos@renyi.hu)}}

\maketitle

\begin{abstract}
			In this work we establish an upper bound on the almost sure convergence rate of 
			a broad class of push-sum algorithms, more precisely for those where message passing events follow a Markovian dynamics, or a Hidden Markov Model. In this context, having a computationally accessible quantity at hand is of top priority; we thus formulate tools to cut down the complexity of obtaining the stated bound.

                        Furthermore, we introduce a transformation to adjust a given push-sum scheme in a way that the bound on the rate provably improves -- and often so does the true convergence rate.

			Theoretical understanding is complemented by numerical results, where the quality of all three directions is demonstrated: for the bound on the rate, for the computational complexity and for the improvement strategy.
		\end{abstract}

\section{Introduction}
\label{sec:intro}

Given the gap between the growth of computational need and the growth of computational capacity of a physical unit, relying on distributed solutions is indispensable. This creates the additional layer of coordination between computational units; a fundamental building block is average consensus aggregating data to their average.

Our current work is on assessing the convergence rate of this element. We focus on the celebrated construction of push-sum, which dates back to the work of \cite{kempe2003gossip} and its follow-up \cite{benezit2010weighted}.
The strength of the protocol is its flexibility and robustness, being able to cope with asynchronous and/or directed communication. Variants for further communication challenges have been designed, e.g.\ for packet drops in \cite{hadjicostis2015robust}.

A crucial question is efficiency, both when the method is used as is or when used as a building block. More precisely, this can be formally captured by the speed of convergence of the actual push-sum construction.

\subsubsection*{Contributions}
In the current work, we focus on three related directions. First, to bound the almost sure (a.s.) exponential convergence rate when Markovian dynamics drive message passing, extending \cite{gerencser2022computable}. Second, to reformulate the bounding quantity to decrease computational complexity exploiting sparse matrix algebra. Third, to propose a method that allows modification of a push-sum scheme to one that often converges faster.

The contributions of the paper can be summarized as follows.

\begin{enumerate}
\item We establish a computable bound on the almost sure exponential convergence rate of push-sum algorithms, where a Hidden Markov Model (HMM) governs the dynamics, i.e.\ a background Markovian process is used to randomize the transmissions. The exact statement is formulated in Theorem \ref{thm:markov} in Section \ref{sec:markovrate}.
\item We support the presented bound with reformulations so that the computational complexity of acquiring it drops significantly. This is presented in Theorem \ref{thm:markovsparse}. Moreover, this reformulation also helps when handling push-sum processes with i.i.d.\ transmissions; this simplified form is stated as Theorem \ref{thm_sp_sparse}; both of these can be found in Section \ref{sec:sparsifying}.
\item We show that a certain ``derandomization'' adjustment to the dynamics can improve its efficiency. More precisely, the convergence rate bounding quantity can provably only decrease, as explained in Theorems \ref{thm:averaged} and \ref{thm:markov_averaged} for i.i.d.\ and HMM cases, respectively. These are stated in Section \ref{sec:derandomize}.
\item In addition, all these theoretical findings are supported by numerical results in Section \ref{sec:numerical} which demonstrate the quality of the bounds and the actual computational effort to acquire them.
\end{enumerate}

\subsubsection*{Related work}

Almost sure convergence has already been proven in a wide range of setups with an exponential decay (even if just implicitly) in \cite{kempe2003gossip}, \cite{benezit2010weighted}, or with communication challenges including packet delay or drop by \cite{hadjicostis2014average}, \cite{gerencser2018push}. Still, the actual rate had not yet been identified at that point.

A key result in this direction is \cite{gerencsr2019tight}, giving a description of the exact a.s.\ rates with the help of the Lyapunov exponents of the associated random matrix product, under very general conditions. This allowed one to put in context and approach the rate with existing tools; however, it is known that this asymptotic quantity is uncomputable in the general case, see \cite{tsitsiklis1997lyapunov}.

The natural continuation in view of these results was to aim for potentially suboptimal, but computable bounds on the rate. Several steps have been taken already. When the communication at each step is randomized in an independent identically distributed (i.i.d.) way, a convincing solution has been found in this direction in \cite{gerencser2022computable}, following and extending the work of \cite{iutzeler2013analysis}, where for a wide range of numerical examples the gap between the computed and simulated rates was an order of magnitude smaller than the rate itself.
An alternative approach has been taken in \cite{gkpushsumsynchr2024}, \cite{gkpushsumpolygen2025} to drastically decrease the computational complexity, leading to similar bounds as before. However, the compromise is that it can be applied only to highly symmetric networks.

It is important to highlight that besides simply averaging out input values, push-sum methods appear as a tool for various applications, recently as the communication layer for distributed optimization, see e.g.\ \cite{nedic2014distributed}, \cite{nedic2017achieving}, \cite{lin2022subgradient} or \cite{bahman_ps_optim2025}.

While push-sum is a popular choice, there are other designs for distributed averaging which can handle asynchronous communication. For example, it is possible to keep linear updates at the price of allowing negative weights \cite{hespanha2018broadcast}. See \cite{hadjicostis2018distributed} for a recent monograph discussing multiple approaches and applications.

\bigskip
The rest of the paper is organized as follows. Later in this section, notation is introduced and a quick overview of the push-sum scheme is presented. Section \ref{sec:markovrate} includes the description of the message passing dynamics considered, the statement and proof of the bound on the convergence rate for the associated push-sum process. In Section \ref{sec:sparsifying} the computational tool is added to help both the current HMM and also the simpler i.i.d.\ case to get the bounding quantity faster. A dynamics adjustment concept and corresponding claim on how it improves can be found in Section \ref{sec:derandomize}. Numerical experiments are presented in Section \ref{sec:numerical} to help intuition and demonstrate the capabilities of the tools built earlier. Conclusions are drawn in Section \ref{sec:conclusions}.

\subsection{Notation, preliminaries}

We will use standard notation: $\mR$ for real numbers and $\mR^N,\mR^{N\times M}$ for vectors and matrices of various sizes.
Probability and expectation will be denoted by $\mP,\mE$ respectively, for the complete probability space we are working with. Various norms will appear, by default we use $\ell_2$ norms for vectors and the induced operator norms for matrices, denoted by $\|\cdot\|$. The Frobenius norm for matrices is denoted by $\|\cdot\|_F$. The Kronecker product of two matrices $A,B$ will be written as $A \otimes B$, with the simplification $A \otimes A = A^{\otimes 2}$ whenever relevant.

Assume a network with $N$ nodes, let $(A_t)_{t\geq 1}$ be a sequence of matrices of size $N\times N$, required to be column-stochastic (thus non-negative) and $x_0 \in \mR^N$ be the initial data. The dynamics required for $(A_t)_{t\geq 1}$ will be described in the next section. The auxiliary weight variable is initialized as $w_0=\mathbf 1 \in \mR^N$. Push-sum algorithms, originally introduced in \cite{kempe2003gossip}, are defined via the following time-inhomogeneous updates:
\begin{align*}
	x_t &= A_t x_{t-1},\\ 
	w_t &=  A_t w_{t-1}.
\end{align*}
The estimate of the average of the initial values at coordinate $i$ is given by the point-wise ratio $z_t^i = x_t^i/w_t^i$.
The goal is to have $z_t^i$ converge to average consensus: $z^i_t\to \bar x:=N^{-1}\mathbf 1^\top x_0$ for all $i$. A.s.\ exponential convergence is known to take place under mild conditions \cite{kempe2003gossip},\cite{gerencsr2019tight}, certifying the usability of the concept.

\section{Rate estimates for Markovian dynamics}
\label{sec:markovrate}
In this section we generalize the tools for bounding the exponential rate of push-sum convergence, relying on \cite{iutzeler2013analysis},\cite{gerencser2022computable}, to the case when the matrices $A_t$ are driven by a Markov process.

More precisely, let $X_t$ be a discrete-time Markov chain on states $[M]=\{1,\ldots,M\}$ and let $P$ denote its transition probability matrix, i.e. 
$p_{mm'} = \mathbb P(X_{t}=m' \mid X_{t-1}=m)$. Suppose that to each state $m$ corresponds a probability measure $\nu_m$ on $ \mR^{N\times N}$.
Assume $A_t$ is obtained via a Hidden Markov Model (HMM), with the dynamics of the joint variable $(X_t, A_t)$ as follows:
\begin{eqnarray}
  \mathbb P(X_{t} = m' \mid X_{t-1} = m) &= &p_{mm'}, \label{eq:mcevolve}\\
  \mathbb P(A_t\in \cdot \mid X_t=m) &=& \nu_m(\cdot),
 \end{eqnarray}
 and naturally $X_t$ and $A_t$ are conditionally independent of the further past. We assume stationary initialization, $X_0$ having law $\pi$, which denotes the stationary distribution of the Markov chain.

 Let us collect the technical conditions needed for the current analysis.
 \begin{assumption}
   \label{ass:markov}
   The matrices $A_t$  are non-negative, allowable and column stochastic. Together they are sequentially primitive, i.e.\ $A_{\psi}\cdots A_1$ is element-wise positive for some random finite $\psi$ a.s. Furthermore, for the minimal non-zero elements $\alpha_t := \min\{A_t^{ij}\mid A_t^{ij} >0\}$ we require $\mE|\log \alpha_1|^{4+\epsilon}<\infty$ for some $\epsilon>0$. The transition matrix $P$ is aperiodic and irreducible. 

\end{assumption}

For the technical analysis we will often need the complete collection of possible transmissions at a specific time, i.e. $(A_{t,1},\ldots,A_{t,M})$, out of which only $A_{t,X_t}$ is the one to be realized.
We are now ready to state the main theorem of this section, to be proved in multiple steps.
 \begin{theorem}
   \label{thm:markov}
   Let Assumption \ref{ass:markov} be in force. With $J=\frac 1N 11^\top$ let
   \begin{align*}
     \bar B &= diag(\mE(A_{1,1}^{\otimes 2}),\ldots,\mE(A_{1,M}^{\otimes 2}))(P^\top\otimes I_{N^2})\\
     \bar B^\circ &= \bar B (I_M\otimes(I-J)^{\otimes 2}).
   \end{align*}
   Then the convergence rate can be almost surely (a.s.) bounded as
   \[
     \limsup_{t\to\infty}\frac 1t \log \sum_{i=1}^N\left|\frac{x_t^i}{w_t^i}-\bar x \right|\leq \frac {\log\rho(\bar B^\circ)}{2}.
   \]
 \end{theorem}

The first main component handles the joint effect of the HMM matrix sequence considered.
 
 \begin{lemma}
   \label{lm:expectedrate}
   Under Assumption \ref{ass:markov}, $\bar B^\circ$ defined as in Theorem \ref{thm:markov}, considering the linear updates in expectation we have
   \[
     \limsup_{t\to \infty} \frac 1t \log \mE\|A_t\cdots A_1(I-J)\| \leq \frac {\log\rho(\bar B^\circ)}{2}.
     \]
 \end{lemma}
 \begin{proof}
   Observe that $A_s(I-J)=(I-J) A_s(I-J)$, thus
   $A_t\cdots A_1(I-J) =  A_t(I-J)\cdots A_1(I-J)$.
   We want to understand the asymptotic rate of
   $
   \|A_t(I-J)\cdots A_1(I-J)\|
   $ in expectation.
   We can reformulate using the Frobenius norm as
   \begin{equation}
     \label{eq:tofrob}
\|A_t(I-J)\cdots A_1(I-J)\|^2\leq \|A_t\cdots A_1(I-J)\|_F^2 = L(A_t^{\otimes 2}(I-J)^{\otimes 2}\cdots A_1^{\otimes 2}(I-J)^{\otimes 2}),     
   \end{equation}
where $L$ is the $\mR^{N^2 \times N^2}\to\mR$ linear map so that $L(A\otimes A)=\|A\|_F^2$. Here $A_t$ follows an HMM structure, we now construct a proxy i.i.d.\ sequence that is easier to handle.

   First, we use that a random function representation of the Markov chain is always possible, see e.g.\ \cite{diaconis1999iterated}, therefore using the (possibly enlarged) underlying probability space $\Omega$ we assume there exists a collection of maps $f_t:[M]\times \Omega \to [M]$ so that
 \[
   X_{t}(\omega) = f_t(X_{t-1}(\omega),\omega)
 \]
 is consistent with the (time-homogeneous) dynamics of \eqref{eq:mcevolve}. We can assume $f_t$ and $X_{t-1}$ to be independent, thus from now on we suppress $\omega$, treating $f_t$ as i.i.d.\ random functions.

Let $B_t$ be composed of $M\times M$ blocks of $N^2\times N^2$ submatrices, all of these blocks being 0 unless $i = f_t(j,\omega)$ for $i,j\in [M]$. For any $i\in [M]$, draw $A_{t,i}$ according to $\nu_i$ independently of $f_t$, and place $A_{t,i}^{\otimes 2}$ in the blocks of row $i$ -- possibly multiple times -- where no 0 has been filled, where some $i = f_t(j,\omega)$ occurs. Technically, we do not need these $M$ matrices to be independent, but we require them to have the same joint law for all $t$. For example, $B_t$ could look like

  \[
   B_t = \left(\begin{array}{ccccc}
     &A_{t,1}^{\otimes 2}&A_{t,1}^{\otimes 2}&&\\
     A_{t,2}^{\otimes 2}&&&&\\
     &&&A_{t,3}^{\otimes 2}&\\
     &&&&\\
     &&&&A_{t,5}^{\otimes 2}
   \end{array}\right).
 \]
By construction, the sequence $(B_t)$ consists of i.i.d.\ variables.
For a more compact notation, let $F_t$ be the binary $M\times M$ matrix representing the mapping $f_t$ when applied from the left. Let $D_t=diag(A_{t,1}^{\otimes 2},\ldots,A_{t,M}^{\otimes 2})$ be a block-diagonal matrix collecting all the available transformation matrices at time $t$. Then we have
\begin{equation}
  \label{eq:Bstruct}
  B_t = D_t (F_t \otimes I_{N^2}).
\end{equation}
 Consider for any vector $v\in \mR^{N^2}$ the augmented vector $\tilde v = e_{X_{t-1}}\otimes v$ using the basis vectors $e_i\in \mR^M$. Thanks to our construction, we get
 \begin{equation}
   \label{eq:Beffect}
   B_t \tilde v = B_t (e_{X_{t-1}} \otimes v) = e_{X_{t}} \otimes (A_{t,X_t}^{\otimes 2}v).
 \end{equation}
 Using this observation we can reformulate the main expression in \eqref{eq:tofrob} as
 \begin{equation}
   \label{eq:AtoB}
   L(A_t^{\otimes 2}(I-J)^{\otimes 2}\cdots A_1^{\otimes 2}(I-J)^{\otimes 2}) =
   (1_M^\top\otimes L)(B_t(I_M\otimes (I-J)^{\otimes 2})\cdots B_1(I_M\otimes (I-J)^{\otimes 2})) \left(e_{X_0}\otimes I^{\otimes 2}\right).
 \end{equation}
 This way we can bound the decay rate of interest in expectation, relying on the fact that $B_t$ are i.i.d.:
 \begin{align*}
   \limsup_{t\to\infty} \frac 1t \log\mE \left\|\prod_{i=1}^tA_i(I-J)\right\| &\leq \limsup_{t\to\infty} \frac 1{2t}  \log \left[ (1_M^\top\otimes L) \left(\prod_{i=1}^t \mE(B_i(I_M\otimes (I-J)^{\otimes 2}))\right)\mE (e_{X_0}\otimes I^{\otimes 2})\right]\\
&\leq \limsup_{t\to\infty} \frac 1{2t}  \log \left[\|1_M^\top\otimes L \|\cdot \|\mE(B_1(I_M\otimes (I-J)^{\otimes 2}))^t\|\cdot \|\mE e_{X_0}\otimes I^{\otimes 2}\|\right]\\ &\leq \frac 12 \log \rho(\mE B_1(I_M\otimes (I-J)^{\otimes 2})).
\end{align*}
Recalling $\mE F_1=P^\top$ the matrix in the final expression coincides with $\bar B^\circ$, allowing us to conclude.
 \end{proof}

  The other main technical component follows, to handle the fluctuations of $w_t$, parallel to \cite[Lemma 4]{gerencser2022computable}:
  \begin{lemma}     Under Assumption \ref{ass:markov} we have 
   \label{lem:w_subexp}
   \begin{equation}
     \label{eq:w_subexp_markov}
     \lim_{t\to \infty} \frac 1t \log \frac 1 {\min_i w_t^i} = 0.     
   \end{equation}
 \end{lemma}
 \begin{proof}
   The statement will be an immediate consequence of \cite[Lemma 44]{gerencsr2019tight} about the similarity of the rows of the matrix products $A_{t:1} := A_tA_{t-1}\cdots A_1$.
   To satisfy the conditions of the cited lemma, the minimal non-zero values $\alpha_n$
   need to be controlled as
   \begin{equation}
     \label{eq:ggmixingcond}
     \mE\left|\sum_{s=1}^t (\log \alpha_s - \mE \log\alpha_s)\right|^q \lesssim t^{q/2},
   \end{equation}
   for all $t$ and some $q>4$. Here $\lesssim$ stands for bounding up to constants, which might depend on the dynamics and $q$, but not on $t$. 
   First, for any finite state ergodic Markov chain a Doeblin condition holds, i.e.\ there exists $k\in\mZ_+$ so that the transition probabilities can be decomposed as $P^k(x,\cdot) = \eps \pi(\cdot)+(1-\eps)\tilde P(x,\cdot)$. Here $\eps\in(0,1)$ and $\tilde P(\cdot,\cdot)$ is again a transition probability matrix with stationary distribution $\pi$.

   For the moment, assume $t=0\mod k$, then we can decompose the left hand side of \eqref{eq:ggmixingcond} as
   \begin{align}
     \mE\left|\sum_{s=1}^t (\log \alpha_s - \mE \log\alpha_s)\right|^q
     &= \mE\left|\sum_{i=1}^k \sum_{s=0}^{t/k-1} (\log \alpha_{sk+i} - \mE \log\alpha_{sk+i})\right|^q \nonumber \\
     &\leq k^{q-1}\sum_{i=1}^k\mE\left|\sum_{s=0}^{t/k-1} (\log \alpha_{sk+i} - \mE \log\alpha_{sk+i})\right|^q \nonumber\\
     &= k^q \mE\left|\sum_{s=0}^{t/k-1} (\log \alpha_{sk+1} - \mE \log\alpha_{sk+1})\right|^q, \label{eq:centeredalphaq}
   \end{align}
   where we applied Jensen's inequality to the convex function $|\cdot|^q$, with the scaling rearranged, then simplified thanks to the stationarity of the process.
   If $t$ had not been divisible by $k$, we would get the analogous expression with a mixture for sums until $\lfloor t/k \rfloor-1$ and $\lceil t/k \rceil-1$. Given that our bounds will be general in the number of summands, it is enough to consider the homogeneous variant above, we will only need a slight increase of the global scaling constant.
   Let us introduce the notation
   \[
     \beta_s = \log \alpha_{sk+1} - \mE \log\alpha_{sk+1}, \qquad \text{for } 0\leq s < \lceil t/k\rceil=:t'.
   \]
   Thanks to the Doeblin condition, we can generate $\beta_s$ by using an $\eps$-coin at each step to decide to reset to an independent stationary variable. This will lead to streaks of dependent variables of length $l_1,l_2,\ldots,l_{r}$, with distribution $Geo(\eps)$, possibly cut for the last to stop at $t'$. Note that the count $r$ is also random. Define the merged terms $\gamma_1,\gamma_2,\ldots$, and their bounding quantities:
   \begin{align*}
     \gamma_i &= \sum_{s=l_1+\ldots+l_{i-1}+1}^{l_1+\ldots+l_{i-1}+l_i}\beta_s,\\
     \tilde\gamma_i &= \sum_{s=l_1+\ldots+l_{i-1}+1}^{l_1+\ldots+l_{i-1}+l_i}|\beta_s|.
   \end{align*}
   Conditioning on $r,l_1,\ldots,l_r$ we have independent, centered $\gamma_i\in L_q$, so by the Marcinkiewicz–Zygmund inequality %XXXCITE
   we have
 \[
     \mE\left(\left|\sum_{i=1}^{r} \gamma_i \right|^q \mid r,l_1,\ldots l_r \right) \lesssim \mE\left(\left(\sum_{i=1}^{r} \gamma_i^2 \right)^{q/2} \mid r,l_1,\ldots l_r \right) \leq \mE\left(\left(\sum_{i=1}^{r} \tilde\gamma_i^2 \right)^{q/2} \mid r,l_1,\ldots l_r \right). 
   \]
   Observe that by the construction of $\tilde\gamma_i$, it can only increase if the number of summands increases.
   So let us extend and increase the expression on the right hand side, starting by completing $l_r$ to follow a proper Geometric distribution (with a slight abuse of notation), then reading the process further until $l_{t'}$ and the corresponding $\tilde\gamma_{t'}$ are defined. As $r\leq t'$ is immediate, we indeed do not decrease the number of summands. In the end, we get the bounding expression below, which again can be simplified using Jensen's inequality.
   \begin{equation}
     \label{eq:conditionalbound}
     \mE\left(\left(\sum_{i=1}^{t'} \tilde\gamma_i^2 \right)^{q/2} \mid l_1,\ldots l_{t'} \right)
     \leq  t'^{q/2-1}\mE\left(\sum_{i=1}^{t'} \tilde\gamma_i^q \mid l_1,\ldots l_{t'} \right)
     =  t'^{q/2-1}\mE\left(\sum_{i=1}^{t'} \tilde\gamma_i^q \mid l_i \right).
   \end{equation}

   For the final term, the summands can be easily bounded using Jensen's inequality and the original moment assumptions:
   \begin{align*}
     \mE(\tilde\gamma_1^{q}\mid l_1) &\leq l_1^{q}\mE(|\beta_1|^{q}) \leq l_1^q2^{q-1}\left(\mE |\log \alpha_1|^q + |\mE \log \alpha_1|^q\right) \lesssim l_1^q.
   \end{align*}
   Substituting into \eqref{eq:conditionalbound} above and taking expectation we get the bound on \eqref{eq:centeredalphaq} as
   \[
      \mE\left|\sum_{s=1}^t (\log \alpha_s - \mE \log\alpha_s)\right|^q \lesssim \mE\left(\sum_{i=1}^{t'}l_i^2\right)^{q/2}\leq {t'}^{q/2}\mE(l_1^{q}) \lesssim t^{q/2},
   \]
   using that all moments of a geometric random variable are finite. Thus we confirm \eqref{eq:ggmixingcond}.
   
   Thanks to the HMM structure, the pattern of non-zeros of $A_{t:1}$ also follows a Hidden Markov chain, with the all-one pattern being the sink. More importantly, for the index of primitivity $\psi_t$ on this finite state space we automatically have $\mE \psi_t <\infty$.
   Then using \cite[Lemma 44]{gerencsr2019tight} we obtain $\lim \frac 1t \log A_{t:1}^{ik}/A_{t:1}^{jk} = 0$. Recalling $w_t = A_{t:1}\mathbf{1}$, the same sub-exponential bound follows for the mixture of the fractions as $\lim \frac 1t \log w_t^{i}/w_t^{j} = 0$. Summing over $i$ yields the claim for any $j$, thus also confirming it for the worst case.
 \end{proof}

 With the above tools at hand, the proof of the theorem above is easy to complete.
 \begin{proof}[Proof of Theorem \ref{thm:markov}]
   We can simplify our target expression as
   \[
  \limsup_{t\to\infty} \frac 1t \log \sum_{i=1}^N\left|\frac{x_t^i}{w_t^i}-\bar x \right|
  \leq \limsup_{t\to\infty}\frac 1t \log \|x_t-w_t\bar x\|_1 - \liminf_{t\to\infty}\frac 1t \log \min_i w_t^i,
\]
and here the second term is 0 thanks to Lemma \ref{lem:w_subexp}. For the first term, using the equivalence of norms, we have
\[
\|x_t-w_t\bar x\|_1 \lesssim  \|x_t-w_t\bar x\| = \|A_t\cdots A_1(x_0-J x_0)\| \leq \| A_t\cdots A_1(I-J)\|\cdot \|x_0\|.
\]
By Lemma \ref{lm:expectedrate}, we have control on the decay rate for the expectation. This implies that for any $l$ the probability $\mP(\sfrac{1}{t} \log \|x_t-w_t\bar x\|_1 >\log\rho(\bar B^\circ)/2+1/l)$ is exponentially decaying in $t$, thus summable, so by the Borel--Cantelli lemma this threshold is passed only finitely many times. Using this for all $l\in\mZ_+$, we get the almost sure rate bound as stated.
 \end{proof}

\section{Computability improvements for rate estimates for i.i.d.\ and Markovian dynamics}
\label{sec:sparsifying}

In this section we provide a computational refinement for the rate estimate of Theorem \ref{thm:markov} together with the preceding work. A generic possibility for bounding the a.s. convergence rate for any i.i.d.\ scheme under mild conditions was initially studied in \cite{iutzeler2013analysis} and strengthened in \cite{gerencser2022computable}. Let us recall the technical assumption and the result on the convergence rate used there.

 \begin{assumption}
   \label{ass:iidsparse}
   The matrices $A_t$ are non-negative, allowable and column stochastic. They form an i.i.d.\ sequence of primitive matrices, i.e.\ an element-wise positive (finite) product can be realized with positive probability. Furthermore for the minimal non-zero elements $\alpha_t := \min\{A_t^{ij}\mid A_t^{ij} >0\}$ we require $\mE|\log \alpha_1|<\infty$.
\end{assumption}

\begin{prop}[\cite{gerencser2022computable}, Theorem 8.]
	\label{prp:Main2022}
        Let $A_t$ be a sequence of matrices so that Assumption \ref{ass:iidsparse} is satisfied.
        Then with $\eta_2 = \log\rho\left(\mE(A_{1}^{\otimes 2})(I-J)^{\otimes 2}\right)$ we have a.s.  
	\begin{align*}
          \limsup_{t\to\infty} \frac 1 t \log  \sum_{i=1}^N\bigg|\frac{x_t^i}{w_t^i} -\bar x\bigg| \le {\frac {\eta_2} 2} < 0.
	\end{align*}
\end{prop}

The cited theorem provides a computable bound for the a.s.\ rate, however via the spectral properties of a matrix of size $N^2\times N^2$ which is possibly filled with non-zero entries as a result of the multiplication by the $(I-J)^{\otimes 2}$ term. We now want to provide a computational shortcut to this useful quantity, exploiting the sparsity $A_1$ possibly provides.
This will be performed by comparing the spectra of the matrices
\[
  \mE(A_{1}^{\otimes 2})(I-J)^{\otimes 2},\mE(A_{1}^{\otimes 2}),\mE(A_{1}).
\]
For example, if vertices have at most $d$ neighbors in the communication graph, then $\mE(A_{1}^{\otimes 2})$ has at most $d^2N^2$ non-zero entries, $\mE(A_{1})$ at most $dN$ of them, indeed demonstrating appealing sparsity.
Once the scheme is discussed for i.i.d.\ processes, we will carry out a similar line of thought for the HMM context of Theorem \ref{thm:markov}.

Let us denote by $\Sp(B)$ the spectrum of a matrix $B$, the multiset of eigenvalues (with algebraic multiplicities), and by $\Sp^\circ(B)$ its restriction to the non-zero values.

\begin{theorem}
  \label{thm_sp_sparse}
  Let Assumption \ref{ass:iidsparse} be in force.
  Then
  \begin{equation}
    \label{eq:spectrum-sparse}
    \mathrm{Sp}^\circ\left(\mE(A_{1}^{\otimes 2})(I-J)^{\otimes 2}\right) = (\mathrm{Sp}^\circ(\mE(A_{1}^{\otimes 2})) + \{1\}) - 2\times \mathrm{Sp}^\circ(\mE(A_{1})),
  \end{equation}
  encoding manipulation of multisets, multiplication corresponding to increasing the multiplicities.
  As a result, the spectral radius can be expressed as
  \begin{equation}
    \label{eq:eta2-spectrum-sparse}
    \eta_2 = \log\rho\left(\mE(A_{1}^{\otimes 2})(I-J)^{\otimes 2}\right) = \log\max\left\{|\lambda| \mid \lambda \in \big(\mathrm{Sp}^\circ(\mE(A_{1}^{\otimes 2})) + \{1\}\big) - 2\times \mathrm{Sp}^\circ(\mE(A_{1}))\right\}.
  \end{equation}
\end{theorem}

\begin{remark}	
By a Cauchy--Schwarz inequality for tensor products, see e.g.\ \cite{gerencser2008lq}, we have
$$
\log\rho(\mE(A_{1})(I-J)) \leq \frac 12\log\rho(\mE(A_{1}^{\otimes 2})(I-J)^{\otimes 2}).
$$
Therefore once $\mathrm{Sp}^\circ(\mE(A_1))$ is at hand, we have a lower bound on $\eta_2$.
Given that $\eta_2$ can be obtained via a difference of multisets,
once the top
\begin{equation}
  \label{eq:topeigcount}
2\cdot\left|\{\lambda \in \mathrm{Sp}^\circ(\mE(A_1)) \mid \rho(\mE(A_{1})(I-J))^2\leq |\lambda|\leq 1\}\right|  
\end{equation}
eigenvalues from $ \mathrm{Sp}^\circ(\mE(A_{1}^{\otimes 2}))$ are computed, a distinct value will appear.
\end{remark}

The theorem and the remark above lead to a numerical algorithm as well.
The following steps are needed schematically. \\

\fbox{\parbox{0.95\linewidth}{
	\textbf{Algorithm 1 (i.i.d.\ case)} 
\begin{enumerate}\setlength\itemsep{1pt}
	\setcounter{enumi}{-1}
\item Choose a numerical error tolerance $\varepsilon>0$.
\item Compute the eigenvalues of the matrix $\mE(A_{1})$.
\item Set $k := 2\cdot |\{\lambda \in \mathrm{Sp}(\mE A_1) \mid \rho\big(\mE A_1(I-J)\big)^2 \leq |\lambda| \leq 1\}|$ and compute the first $k$ eigenvalues of largest magnitude of the possibly sparse matrix $\mE(A_1^{\otimes2})$. 
\item Attempt to match $2\times \mathrm{Sp}(\mE(A_1))$ with $ \mathrm{Sp}(\mE(A_{1}^{\otimes 2})) \cup \{1\}$ in decreasing order of magnitudes.
\item Let $\mu_1\in \mathrm{Sp}(\mE(A_1)), \mu_2\in  \mathrm{Sp}(\mE(A_{1}^{\otimes 2}))$ be the first occasion when the difference is larger than $\varepsilon$.
\item Return $\eta_2 = \log|\mu_2|$.
\end{enumerate}}} \\

\begin{proof}[Proof of Theorem \ref{thm_sp_sparse}]
  Within $\mR^{N^2} = \mR^{N}\otimes \mathbb R^N$, consider the subspace $U = (1 \otimes \mR^N) + (\mR^N \otimes 1)$, of dimension $2N-1$. We will analyze the action of $\mE(A_{1}^{\otimes 2})$ in the block matrix structure of the decomposition $\mR^{N^2} = U \oplus U^\perp$.
  Consider any $z\in \mR^N$, and test the action on $1^\top\otimes z^\top\in U$, relying on the matrices being column-stochastic:
 \begin{equation}
   \label{eq:colstoch-exploit}
   1^\top\otimes z^\top\mE(A_{1}^{\otimes 2}) =   \mE(1^\top A_{1}\otimes z^\top A_{1})
  =   1^\top\otimes z^\top \mE(A_{1}) \in U.
\end{equation}
The analogous effect holds for vectors of the form $z^\top\otimes 1^\top\in U$, thus we see $U$ is an invariant subspace w.r.t.\ $\mE(A_{1}^{\otimes 2})$ (applied from the right). Let us write $\mE A_1^{\otimes 2}$  in a block matrix form, with blocks corresponding to projections onto the $U$ or $U^\perp$ space, i.e. if $\mathcal P_U$ denotes the orthogonal projection onto the subspace $U$, then 
$\mathcal A_{11}= \mathcal P_U \mE (A_1^{\otimes2}) \mathcal P_U$, $\mathcal A_{12} = \mathcal P_U \mE(A_1^{\otimes2}) (I-\mathcal P_U)$, $\mathcal A_{21} =(I-\mathcal P_{U})  \mE (A_1^{\otimes 2})\mathcal P_U  $ and $\mathcal A_{22} =(I-\mathcal P_U) \mE(A_1^{\otimes2}) (I-\mathcal P_U)  $, so

\begin{equation*}
	\mE A_1^{\otimes 2} = \begin{bmatrix}
		\mathcal A_{11} & 0 \\
		\mathcal A_{21} & \mathcal A_{22}
	\end{bmatrix} 
\end{equation*}
since $U$ is an invariant subspace, implying that $\mathcal A_{12} = 0$. Due to its block-triangular shape, the spectrum of $\mE A_1^{\otimes2}$ consists of values from the spectra of $\mathcal A_{11}$ and $\mathcal A_{22}$. 
From \eqref{eq:colstoch-exploit} again, $\mathrm{Sp}(\mathcal A_{11}) = 2\times \mathrm{Sp}(\mE(A_1)) - \{1\}$, as we can use it with vectors of the form $1^\top\otimes z^\top$ and also $z^\top\otimes 1^\top$, for all $z\in\mR^N$, but in this way double-counting $1^\top\otimes 1^\top$ with eigenvalue 1.
Now let us view $(I-J)^{\otimes 2}$ in this block structure. Here $(1^\top\otimes z^\top)(I-J)^{\otimes 2} = 0$ for any $z\in\mR^N$, similarly for $z^\top\otimes 1^\top$, also when applied from the right. However, for any $y,z\perp 1$, we have $(y^\top\otimes z^\top)(I-J)^{\otimes 2} = y^\top\otimes z^\top$, whence
\[
	(I-J)^{\otimes2} = \begin{bmatrix}
		0 & 0 \\
		0&  I_{N^2-2N+1}
	\end{bmatrix}.
\]
This immediately leads to
$$ \mE(A_1^{\otimes2}) (I-J)^{\otimes2} = \begin{bmatrix}
	0 & 0 \\ 0 & \mathcal A_{22}
\end{bmatrix}, $$
thus $\mathrm{Sp}^\circ\bigg(\mE A_1^{\otimes2} (I-J)^{\otimes2}\bigg)  = \mathrm{Sp}^\circ(\mathcal A_{22}) = \mathrm{Sp}(\mE A_1^{\otimes2}) - \mathrm{Sp}(\mathcal A_{11}) = \mathrm{Sp}(\mE A_1^{\otimes2}) - \big(2\times \mathrm{Sp}(\mE A_1) - \{1\} \big)$.

\end{proof}

Let us now turn to the case of Markovian transmissions, where similar phenomena can be observed, but it becomes slightly more complicated.

 \begin{theorem}
   \label{thm:markovsparse}
   In the setting of Theorem \ref{thm:markov}, let us extend the list of matrices considered as
   \begin{align*}
     \bar B^\circ &= diag(\mE(A_{1,1}^{\otimes 2}),\ldots,\mE(A_{1,M}^{\otimes 2}))(P^\top\otimes (I-J)^{\otimes 2}),\\
     \bar B &= diag(\mE(A_{1,1}^{\otimes 2}),\ldots,\mE(A_{1,M}^{\otimes 2}))(P^\top\otimes I_{N^2}),\\
     \bar C^\circ &= diag(\mE A_{1,1},\ldots,\mE A_{1,M})\cdot (P^\top\otimes (I-J)).\\
     \bar C &= diag(\mE A_{1,1},\ldots,\mE A_{1,M})\cdot (P^\top\otimes I).
   \end{align*}    
   Then we can relate the spectra as
   \begin{align*}
     \Sp^\circ(\bar B^\circ) &= (\Sp^\circ(\bar B) + \Sp^\circ(P)) - 2\times(\Sp^\circ(\bar C)),\\
     \Sp^\circ(\bar C^\circ) &= \Sp^\circ(\bar C) - \Sp^\circ(P).
   \end{align*}
 \end{theorem}
 Once again we can hope for sparse matrices to appear. If all vertices have at most $d$ neighbors in the communication graph, and all states $m$ have non-zero transition probabilities to at most $D$ other states, then it is easy to see that $\bar B,\bar C,P$ have at most $d^2DN^2M,dDNM,DM$ non-zero entries, respectively.

 \begin{proof}[Proof of Theorem \ref{thm:markovsparse}]
   The argument is closely related to the proof of Theorem \ref{thm_sp_sparse} before, here we emphasize the additional layer of complexity. We again rely on identifying invariant left subspaces. Let us define 
   \begin{align*}
     U &= \langle\mR^M\otimes\mR^N\otimes 1\rangle \subset \mR^{N^2M},\\
     U' &= \langle\mR^M\otimes 1\otimes\mR^N\rangle \subset \mR^{N^2M}.
   \end{align*}
   For any vector $z\in \mR^{MN}$ we have $z\otimes 1\in U$, and we can check the action of $\bar B$, relying on $A_{1,j}$ being column-stochastic:
   \begin{equation}
     \label{eq:uinvariant}
   \begin{aligned}
     (z^\top\otimes 1^\top) &diag(\mE(A_{1,1}^{\otimes 2}),\ldots,\mE(A_{1,M}^{\otimes 2}))(P^\top\otimes I_{N^2})\\
     &= (z^\top diag(\mE A_{1,1},\ldots,\mE A_{1,M})\otimes 1^\top)(P^\top\otimes I_{N^2})\\
     &= (z^\top diag(\mE A_{1,1},\ldots,\mE A_{1,M})(P^\top\otimes I))\otimes 1^\top \in U,
   \end{aligned}  
   \end{equation}
   thus we see $U$ is indeed an invariant subspace. The same reasoning holds for $U'$ by the symmetry of the roles of the two $\mR^N$ components.

   It follows that the intersection $U\cap U'= \mR^M\otimes 1_{N^2}$ is also an invariant subspace. For $z\otimes 1_{N^2}\in U\cap U'$ the action we see is
   \begin{equation}
     \label{eq:uuinvariant}
     \begin{aligned}
     (z^\top\otimes 1_{N^2}^\top) diag(\mE(A_{1,1}^{\otimes 2}),\ldots,\mE(A_{1,M}^{\otimes 2}))(P^\top\otimes I_{N^2}) =
     (z^\top\otimes 1_{N^2}^\top) (P^\top\otimes I_{N^2}) =  (z^\top P^\top)\otimes 1_{N^2}^\top.
   \end{aligned}
   \end{equation}
   
   We can write $\bar B$ in a block matrix form by changing coordinates according to the subspaces just identified: $U\cap U', U^\perp\cap U', U\cap U'^\perp, (U+U')^\perp$, which has the following pattern:
\begin{equation*}
	\bar B= \begin{bmatrix}
		\mathcal B_{11} & 0 & 0 & 0\\
                  \mathcal B_{21} & \mathcal B_{22} & 0 & 0\\
                  \mathcal B_{31} & 0 & \mathcal B_{33} & 0 \\
                  \mathcal B_{41} & \mathcal B_{42} &  \mathcal B_{43} &  \mathcal B_{44}
	\end{bmatrix} .
\end{equation*}

The matrix of interest $\bar B^\circ$ is obtained by further multiplying by $I_M\otimes (I-J)^{\otimes 2}$. Using the coordinates above, this can be written as
\begin{equation*}
  \bar B^\circ = \bar B (I\otimes (I-J)^{\otimes 2})= \begin{bmatrix}
		\mathcal B_{11} & 0 & 0 & 0\\
                  \mathcal B_{21} & \mathcal B_{22} & 0 & 0\\
                  \mathcal B_{31} & 0 & \mathcal B_{33} & 0 \\
                  \mathcal B_{41} & \mathcal B_{42} &  \mathcal B_{43} &  \mathcal B_{44}
                                                  \end{bmatrix} \cdot
                                                  \begin{bmatrix}
                                                    0 & 0 & 0 & 0\\
                                                    0 & 0 & 0 & 0\\
                                                    0 & 0 & 0 & 0\\
                                                    0 & 0 & 0 & I_{M(N-1)^2}
                                                  \end{bmatrix}
                                                  = \begin{bmatrix}
                                                    0 & 0 & 0 & 0\\
                                                    0 & 0 & 0 & 0\\
                                                    0 & 0 & 0 & 0\\
                                                    0 & 0 & 0 & \mathcal B_{44}
	\end{bmatrix}.
\end{equation*}
As the matrices are block-triangular, we immediately get that the difference in the spectrum is replacing $\Sp(\mathcal B_{11})+\Sp(\mathcal B_{22})+\Sp(\mathcal B_{33})$ with zeros.

In the same coordinate base, restricting to the invariant subspace $U$, we see that the spectrum there is $\Sp(\mathcal B_{11})+\Sp(\mathcal B_{33})$. At the same time, by \eqref{eq:uinvariant} we see that the action corresponds to that of $\bar C$. A similar correspondence holds on $U'$. Furthermore, based on \eqref{eq:uuinvariant}, we observe the action of $P^\top$ on $U\cap U'$. Thus for the spectra we get
\begin{align*}
  \Sp(\bar C) &= \Sp(\mathcal B_{11})+\Sp(\mathcal B_{22})\\
   &= \Sp(\mathcal B_{11})+\Sp(\mathcal B_{33}),\\
  \Sp(P) &= \Sp(\mathcal B_{11}).
\end{align*}

This allows us to directly express the spectra: we need to add $\Sp(P)$ and remove $\Sp(\bar C)$ twice to cancel $\Sp(\mathcal B_{11})+\Sp(\mathcal B_{22})+\Sp(\mathcal B_{33})$, confirming the stated expression for $\Sp^\circ(\bar B^\circ)$.

For the claim on $\Sp(\bar C^\circ)$, a similar but simpler reasoning is sufficient. The invariant left subspace of $\bar C$ to be used is $\mR^M\otimes 1$, for any $z\in \mR^M$ the action on $z\otimes 1$ is
\[
  (z^\top\otimes 1^\top) diag(\mE A_{1,1},\ldots,\mE A_{1,M})\cdot (P^\top\otimes I) =
  (z^\top\otimes 1^\top) (P^\top\otimes I) =   (z^\top P^\top)\otimes 1^\top.
\]
Now this is the subspace mapped to 0 by the additional $I_M\otimes (I-J)$ factor of $\bar C^\circ$, thus by the equation above  $\Sp(P)$ needs to be removed from the spectrum, confirming the claim stated for $\Sp^\circ(\bar C^\circ)$.
   
 \end{proof}

As a pleasant consequence we see that the convergence rate bound is indeed meaningful:
 
 \begin{corollary}
   In the setting of Theorem \ref{thm:markov} for the bounding quantity therein we have $\rho(\bar B^\circ)<1$.
 \end{corollary}

 \begin{proof}
   As $\rho(\bar B^\circ)\leq 1$, we simply need to check that the multiplicity of eigenvalue 1 is 0, and that there is no other eigenvalue of modulus 1. All the matrices $\bar B,\bar C,P$ are non-negative, column-stochastic. Once we confirm all three are irreducible and aperiodic we are done, as by Theorem \ref{thm:markovsparse} the multiplicity of eigenvalue 1 for $\bar B^\circ$ is $1+1-2=0$, and such matrices cannot have other eigenvalues on the unit circle.

   For $P$, irreducibility and aperiodicity are by assumption on the Markov chain. In particular $P^{T_P}$ is element-wise positive for large enough constant time $T_P$.
   
   For $\bar C$, consider the i.i.d.\ random matrix sequence $C_t$ similarly to $B_t$ used in the proof of Lemma~\ref{lm:expectedrate}:
   \[
     C_t = diag(A_{t,1},\ldots,A_{t,M})(F_t \otimes I_N).
   \]
   These matrices have $\mE C_t = \bar C$.

   Recall that by sequential primitivity there is an almost surely finite (random) time $\psi$ so that $A_\psi\cdots A_1$ is element-wise positive. There is thus a large enough $T_\psi$ so that $\mP(\psi\leq T_\psi \mid X_0=j)>0$ for all $j\in [M]$.

   We will see that the $(i,j)$ block of $C_{T_P+T_\psi}\cdots C_1$ is element-wise positive with positive probability. Let $X_0=j$, then by the comment above $A_{T_\psi}\cdots A_1$ is element-wise positive with positive probability, with a (random) terminal state $X_{T_\psi}$ for the Markov chain. Thanks to the aperiodicity note on $P$ above, there is a positive probability route from $X_{T_\psi}$ to $i$ in $T_P$ steps for the driving Markov chain. Importantly, thanks to $A_t$ being allowable, $A_{T_P+T_\psi}\cdots A_1$ is still element-wise positive on this event, which is now the $(i,j)$ block of the product of interest.

   Now that all the blocks of $C_{T_P+T_\psi}\cdots C_1$ are element-wise positive with positive probability, it follows that $\mE(C_{T_P+T_\psi}\cdots C_1) = \bar C ^{T_P+T_\psi}$ is an element-wise positive matrix, confirming the irreducibility and aperiodicity of $\bar C$.

   For $\bar B$, the same proof works as for $\bar C$ above. We get element-wise positive $(A_{T_\psi}\cdots A_1)^{\otimes 2}$ with positive probability for any $X_0=j$, so that $B_{T_P+T_\psi}\cdots B_1$ has a fully positive block (of size $N^2\times N^2$) at any position $(i,j)$ with positive probability, confirming irreducibility and aperiodicity as above.
 \end{proof}

Our goal is once again to phrase an algorithmic way to compute the convergence rate bound. One more component is to capture a lower bound on this quantity, as follows.
\begin{lemma}
\label{lm:BC_rhos}
  With the notations of Theorem \ref{thm:markov} and Theorem \ref{thm:markovsparse} we have $\rho(\bar C^\circ)^2\leq \rho(\bar B^\circ)$.
\end{lemma}
The proof is provided in Appendix \ref{sec:tensorcs}.

Thanks to Theorem \ref{thm:markovsparse} we can bound the spectral radius of $\bar B^\circ$ via the spectrum of $\bar C, P$, which are sparse matrices in a number of cases, as indicated above. Moreover, the number of eigenvalues necessary to compute from the largest (sparse) matrix $\bar B$ can be bounded by $k_B$ defined as
\begin{align*}
  k_C &:= |\{\lambda \in \mathrm{Sp}(\bar C) \mid \rho\big(\bar C^\circ\big)^2 \leq |\lambda| \leq 1\}|,\\
  k_P &:= |\{\lambda \in \mathrm{Sp}(P) \mid \rho\big(\bar C^\circ\big)^2 \leq |\lambda| \leq 1\}|,\\
  k_B &:= 2k_C-k_P+1.
\end{align*}
These quantities immediately follow from the relation of spectra in Theorem \ref{thm:markovsparse} and the threshold identified in Lemma \ref{lm:BC_rhos}. With these at hand, we can compose the following schematic algorithm.\\

\fbox{\parbox{0.95\linewidth}{
	\textbf{Algorithm 2a (HMM case)} 
\begin{enumerate}\setlength\itemsep{1pt}
	\setcounter{enumi}{-1}
\item Choose a numerical error tolerance $\varepsilon>0$.
\item Compute the eigenvalues of matrices $\bar C,P$.
\item Attempt to match $\mathrm{Sp}(\bar C)$ with $\mathrm{Sp}(P)$ in decreasing order of magnitudes.
  \item Let $\alpha\in\mathrm{Sp}(\bar C)$ be the first not matched within tolerance $\varepsilon$.
\item Set $k_C := |\{\lambda \in \mathrm{Sp}(\bar C) \mid |\alpha|^2 \leq |\lambda| \leq 1\}|$ and
  $k_P := |\{\lambda \in \mathrm{Sp}( P) \mid |\alpha|^2 \leq |\lambda| \leq 1\}|$.
\item With $k_B:= 2k_C-k_P+1$ compute the first $k_B$ eigenvalues of largest magnitude of $\bar B$. 
\item Attempt to match $2\times \mathrm{Sp}(\bar C)$ with $\mathrm{Sp}(\bar B) + \mathrm{Sp}(P)$ in decreasing order of magnitudes.
\item Let $\mu\in \mathrm{Sp}(\bar B)$ be the first occasion when there is no match within $\varepsilon$.
\item Return $\log|\mu|$. \end{enumerate}}} \\ 

An alternative online variant can be formulated, requiring a \emph{stream} of eigenvalues supplied in decreasing magnitude rather than a collection at once.\\

\fbox{\parbox{0.95\linewidth}{
	\textbf{Algorithm 2b (HMM case, iterative)} 
\begin{enumerate}\setlength\itemsep{1pt}
	\setcounter{enumi}{-1}
      \item Choose a numerical error tolerance $\varepsilon>0$.
      \item Initialize $\lambda_B=\lambda_C=\lambda_P=\{1\}$.
      \item Set $\mu = \min{\{ |\lambda | \mid \lambda\in\lambda_B\}}$.
      \item If $\mu > \min{\{ |\lambda | \mid \lambda\in\lambda_C\cup\lambda_P\}} +\varepsilon$:
      return $\log\mu$.
      \item Compute all eigenvalues of magnitude $\mu$ for $\bar B,\bar C,P$, include them in $\lambda_B,\lambda_C,\lambda_P$.
      \item Having all is certified by finding the first smaller one for all three: $\mu^-_B,\mu^-_C,\mu^-_P$.
      \item Create minimal distance matching between $2\times\lambda_C$ and $\lambda_B+\lambda_P$.
      \item If the minimal distance is $>\varepsilon$:
      return $\log \mu$.
    \item Reset to $\lambda_B=\{\mu^-_B\},~\lambda_C=\{\mu^-_C\},~\lambda_P=\{\mu^-_P\}$.
      \item Repeat from Step 2.
\end{enumerate}}}

\section{Improvement of dynamics using derandomization}
\label{sec:derandomize}
Here we establish that reducing randomness in the transmission protocol is likely to improve its efficiency. More precisely, we show the decrease of the bounding quantities studied in Proposition \ref{prp:Main2022} and Theorem \ref{thm:markov}.

First consider the case of i.i.d.\ steps. Let $\cG_t$ be a sequence of $\sigma$-algebras concerning transmissions at time $t$. That is, we require $\sigma(A_1,\cG_1),\sigma(A_2,\cG_2),\ldots$ to be independent. We also require them to represent the same information for all $t$, in particular that $\mE(A_t \mid \cG_t)$ are identically distributed for all $t$. For convenience, we introduce the notation
\[
    \AG{t} := \mE(A_t \mid \cG_t).
\]

\begin{lemma}
  \label{lm:cone}
  Define the cone $\cC$ via its generators as $\cC := \langle x^{\otimes 2}, x\in \mR^N\rangle \subset \mR^{N^2}$. Then $\cC$ is an invariant cone for any $X\in \langle A^{\otimes 2}, A \in \mR^{N\times N}\rangle$. Also, for any $a,b\in\cC$ there holds $\|a+b\|_2\geq \|a\|_2$.
\end{lemma}
\begin{proof}
  For any $A,x$, we have $A^{\otimes 2} x^{\otimes 2} = (Ax)^{\otimes 2}\in \cC$. Taking a non-negative mixture of possible $A$ and possible $x$ creates a mixture of the product but remains in $\cC$ as it is a cone, confirming the first statement.

  Considering the second claim, for any $u,v\in\mR^N$,
  \[
    \langle u^{\otimes 2},v^{\otimes 2}\rangle =    \langle u,v\rangle^2 \geq 0,
  \]
  and similarly for non-negative mixtures, so for any $a,b\in\cC$ we have $\langle a,b\rangle\geq 0$.
  For the differences of norms, we can then write
  \[
    \|a+b\|^2-\|a\|^2 = 2\langle a,b\rangle + \|b\|^2\geq 0.
  \]
\end{proof}

Define $U=\cC-\cC$, the subspace of vectors symmetric in their two components. A straightforward consequence of the Lemma above is that $U$ is an invariant subspace for matrices in $\langle A^{\otimes 2}, A \in \mR^{N\times N}\rangle$.

\begin{lemma}
  \label{lm:radius_restricted}
  Let $Y$ be a random matrix in $\mR^{N\times N}$. Denote by $\|\cdot\|_{2,U}$ and $\rho_U(\cdot)$ the induced $L_2$ norm and the spectral radius when the operator is restricted to $U$, respectively. Then there holds
  \begin{equation*}
    \|\mE(Y^{\otimes 2})\|_2 \lesssim \|\mE(Y^{\otimes 2})\|_{2,U}\lesssim \|\mE(Y^{\otimes 2})\|_2.
  \end{equation*}
  Additionally,
  \begin{equation*}
    \rho_U(\mE(Y^{\otimes 2}))=\rho(\mE(Y^{\otimes 2})).
  \end{equation*}
\end{lemma}
\begin{proof}
  The induced norm on $U$ is bounded above by the induced norm on $\mR^{N^2}$ as it is the same operator restricted, the non-trivial step is the other direction.

  Consider the test vector $u = \sum_{i=1}^Ne_i^{\otimes 2}\in U$. As $\|u\|_2^2=N$, we can bound $\|\mE(Y^{\otimes 2})\|_{2,U}$, the operator norm when acting on $U$ as
  \begin{align*}
    N\|\mE(Y^{\otimes 2})\|^2_{2,U}\geq u^\top\mE(Y^{\top\otimes 2})\mE(Y^{\otimes 2})u=\mE(u^\top (Y^\top\tilde Y)^{\otimes 2} u),
  \end{align*}
  where $\tilde Y$ is an independent copy of $Y$. Using the definition of $u,\tilde Y$, the cyclic invariance and linearity of the trace operator this can be rewritten as
  \begin{align*}
    \mE\left(\sum_{i,j=1}^N(Y^\top\tilde Y)^2_{i,j}\right)
    &=\mE (\tr Y^\top\tilde Y \tilde Y^\top Y)
    = \mE (\tr \tilde Y \tilde Y^\top YY^\top)
= \tr\mE ( \tilde Y \tilde Y^\top YY^\top)
= \tr\mE (\tilde Y\tilde Y^\top) \mE(YY^\top)\\
    &= \tr (\mE(YY^\top)^2).
  \end{align*}

  For comparison let us consider the Frobenius norm of $\mE(Y^{\otimes 2})$, which is sufficient thanks to the equivalence of norms.
  Using the same observations and a Cauchy-Schwarz inequality, as the trace realizes an inner product, we obtain
  \begin{align*}
    \|\mE(Y^{\otimes 2})\|^2_F = \mE(\tr  \tilde Y^{\otimes 2}Y^{\top \otimes 2}) = \mE((\tr  \tilde Y Y^\top)^2)
    \leq \mE((\tr  \tilde Y \tilde Y^\top)(\tr  Y Y^\top))=\mE(\tr  Y Y^\top)^2=(\tr\mE(Y Y^\top))^2.
  \end{align*}

  At the same time, we have  $N\tr (\mE(YY^\top)^2)\geq (\tr\mE(Y Y^\top))^2$ by comparing the $\ell_2$ and $\ell_1$ norms of the spectra of the PSD matrix $\mE(YY^\top)$. Combining with the above we arrive at
  \begin{equation*}
     \frac{N^2\|\mE(Y^{\otimes 2})\|^2_{2,U}}{\|\mE(Y^{\otimes 2})\|^2_F}\geq 1,
   \end{equation*}
   confirming the claim relating the norms.
   
   For the spectral radii, take the product of $k$ independent copies of $Y$ and apply the claim on norms for it. Then the matrix studied becomes $\mE(Y^{\otimes 2})^k$, for which the above calculation and comparison of norms also holds. If there had been a gap between the spectral radii, then the ratio of norms would asymptotically vanish, contradicting the bound obtained, thus the spectral radii must be equal.
\end{proof}

\begin{theorem}
  \label{thm:averaged}
  Let Assumption \ref{ass:iidsparse} be in force. For the convergence rate bounds of Proposition \ref{prp:Main2022} for the two (i.i.d.) matrix sequences we have
  \[
    \rho\left(\mE(\AG{1}^{\otimes 2})(I-J)^{\otimes 2}\right) \leq     \rho\left(\mE(A_{1}^{\otimes 2})(I-J)^{\otimes 2}\right).
  \]
  Consequently, if Assumption \ref{ass:iidsparse} also holds for $\AG{t}$, the corresponding modified algorithm has a better a.s.\ convergence rate bound via Proposition \ref{prp:Main2022}.
\end{theorem}
\begin{proof}
  Let us look at the following term comparing the two random matrix distributions:
  \begin{equation}
    \label{eq:avgdiff}
  \begin{aligned}
    \Delta &:= \mE((A_1-\mE(A_1\mid \cG_1))^{\otimes 2})(I-J)^{\otimes 2}\\
      &= \mE(\mE(A_1\otimes A_1 -A_1\otimes\mE(A_1\mid \cG_1) - \mE(A_1\mid \cG_1)\otimes A_1 + \mE(A_1\mid \cG_1)\otimes\mE(A_1\mid \cG_1) \mid \cG_1))(I-J)^{\otimes 2} \\
      &= \mE(A_1^{\otimes 2}-\mE(A_1\mid \cG_1)^{\otimes 2})(I-J)^{\otimes 2}\\
      &= \mE(A_1^{\otimes 2}(I-J)^{\otimes 2}) - \mE(\AG{1}^{\otimes 2}(I-J)^{\otimes 2}).
  \end{aligned}    
  \end{equation}
  Here we rely on $\mE(A_1\otimes\mE(A_1\mid \cG_1) \mid \cG_1) = \mE(A_1\mid \cG_1)\otimes \mE(A_1\mid \cG_1)$, thanks to the $\cG_1$-measurability of the second term.
  From the above, $\Delta$ expresses the difference of the matrices needed for accessing the spectral radii of interest.

  When looking for the vectors that determine the spectral radius, it will be sufficient to search for symmetric vectors in $U$ thanks to Lemma \ref{lm:radius_restricted}.

  By Lemma \ref{lm:cone} above $\cC\subset U$ is an invariant cone for $\mE(\AG{1}^{\otimes 2}(I-J)^{\otimes 2})$, and thanks to $\cC$ being a proper pointed cone within $U$, the Krein--Rutman theorem is applicable.
  That is, there exists a $v\in\cC$ principal (right) eigenvector which realizes the spectral radius. Thanks to the decomposition above, we have
  \[
    \mE(A_1^{\otimes 2}(I-J)^{\otimes 2})v = \mE(\AG{1}^{\otimes 2}(I-J)^{\otimes 2})v + \Delta v.
  \]
  By the definition of $\Delta$, again $\cC$ is invariant for it, thus iterating we get
  \begin{equation}
    \label{eq:spr_iterate}
    \left(\mE(A_1^{\otimes 2}(I-J)^{\otimes 2})\right)^kv = \left(\mE(\AG{1}^{\otimes 2}(I-J)^{\otimes 2}) + \Delta\right)^kv = \left(\mE(\AG{1}^{\otimes 2}(I-J)^{\otimes 2})\right)^kv + z_k,
  \end{equation}
  with some $z_k\in \cC$. This additional term can only increase the norm as shown in Lemma \ref{lm:cone}.
  Thus in \eqref{eq:spr_iterate} the norm can only increase from the expression on $\AG{1}$ to the one on $A_1$, confirming the claim on the comparison of spectral radii.
\end{proof}

After having shown the abstract comparison, let us present a few illustrative examples.

\begin{example}
  Given a connectivity graph $G$, at every step $t$ a random vertex $v_t$ sends a packet of random proportion $q_t$ to a randomly chosen neighbor $u_t$, which is described by $A_t$. Assume $q_t$ is independent of the rest.
 By setting $\cG_t=\sigma(v_t,u_t)$, we obtain $\AG{t}$ where the same communication events happen but with fixed proportion $\mE q_t$.
\end{example}

Observe that in the example below $\cG_t\nsubseteq \sigma(A_t)$, demonstrating the flexibility of the derandomization approach.
\begin{example}
  In a similar setting as before, let us assume a fixed message proportion $q$. The recipient can be chosen in a redundant way, by two random neighbors $u^{(1)}_t,u^{(2)}_t$ and $\xi_t$ uniform on $\{1,2\}$ to decide which one to pick. By setting $\cG_t=\sigma(v_t,u^{(1)}_t,u^{(2)}_t)$, we obtain $\AG{t}$ corresponding to a process where at every step the random vertex $v_t$ sends two packets of proportion $q/2$ each to two neighbors $(u^{(1)}_t,u^{(2)}_t)$.
\end{example}

We now turn to the analogous question when the transmissions follow Markovian dynamics. The fundamental logic is very similar to the one above, but slightly more technical care is needed.

In this context, let $\cG_t$ be a sequence of $\sigma$-algebras concerning transmissions at time $t$, i.e.\ we require that $\sigma(A_{t,1},\ldots,A_{t,M},\cG_t)$ are independent for all $t$. The $\sigma$-algebras should 
represent the same structure at every $t$, we assume that $(\mE(A_{t,1}\mid \cG_t),\ldots,\mE(A_{t,M}\mid \cG_t))$ are identically distributed (as $t$ varies). Moreover, $\cG_t$ should not be related to the Markovian evolution, thus $\cG_t$ and $F_t$ are required to be independent, recalling that $F_t$ is the matrix of state changes at time $t$, see \eqref{eq:Bstruct}.
Note that for these conditions we use all random $A_{t,i}$ matrices, not only the ones actually applied.
In the simplest setting, $\cG_t = \sigma(\cG_{t,1},\ldots,\cG_{t,M})$, with some time-homogeneous $\cG_{t,m}\subseteq \sigma(A_{t,m})$.

Define our matrices of interest for bounding the rate for the conditional expectations:
\begin{equation}
  \label{eq:markov-nonrandom}
  \begin{aligned}
  B_t^\cG &= diag(\mE(A_{t,1}\mid \cG_t)^{\otimes 2},\ldots,\mE(A_{t,M}\mid \cG_t)^{\otimes 2})(F_t\otimes I_{N^2}),\\
  {\bar B}^{\circ\cG} &= \mE B^\cG_1(I_M\otimes (I-J)^{\otimes 2}).
\end{aligned}
\end{equation}

\begin{theorem}
  \label{thm:markov_averaged}
  Let Assumption \ref{ass:markov} be in force. For the matrices defined in Theorem \ref{thm:markov} and \eqref{eq:markov-nonrandom} we have
  \[
\rho({\bar B}^{\circ\cG})\leq \rho({\bar B}^{\circ}).
\]
Consequently, if Assumption \ref{ass:markov} also holds for the process driven by $\mE(A_{t}\mid \cG_t)$,
derandomization again provides an algorithm with an improved bound on the a.s.\ convergence rate.
\end{theorem}
\begin{proof}
  Exactly as in \eqref{eq:avgdiff}, for every $i\in [M]$ we can define
\[
  \Delta_i := \mE((A_{1,i}-\mE(A_{1,i}\mid \cG_1))^{\otimes 2})(I-J)^{\otimes 2}
= \mE(A_{1,i}^{\otimes 2}(I-J)^{\otimes 2}) - \mE(\AG{{1,i}}^{\otimes 2}(I-J)^{\otimes 2}).
\]
For the main derandomized matrix of interest we can write
\begin{align*}
  {\bar B}^{\circ\cG} &= \mE\left[diag(\mE(A_{1,1}\mid \cG_1)^{\otimes 2},\ldots,\mE(A_{1,M}\mid \cG_1)^{\otimes 2})(F_1\otimes (I-J)^{\otimes 2})\right]\\
                   &=diag(\mE[\mE(A_{1,1}\mid \cG_1)^{\otimes 2}(I-J)^{\otimes 2}],\ldots,\mE[\mE(A_{1,M}\mid \cG_1)^{\otimes 2}(I-J)^{\otimes 2}])(P^\top\otimes I_{N^2}).
\end{align*}
Here we use the independence of $\cG_t,F_t$ and linearity to rearrange the expectation and the $(I-J)^{\otimes 2}$ terms.
We can now combine with the definition of $\bar B^\circ$ to express the difference of the matrices of interest:
\[
  \bar B^\circ-{\bar B}^{\circ\cG}
  = diag(\Delta_1,\ldots,\Delta_M) (P^\top\otimes I_{N^2}).
\]

A second look at Lemma \ref{lm:radius_restricted} will reveal that we can restrict our attention to $U^M$ for $\bar B^\circ$ (and similarly for ${\bar B}^{\circ\cG}$). Let us compare the induced norms of $(\bar B^\circ)^k$ on $U^M$ and on $\mR^{N^2M}$.
Using the test vector $1_M\otimes u\in U^M$ we have
\[
  \|(\bar B^\circ)^k(1_M\otimes u)\|_2\geq \max_{i,j\in [M]} \|(e_i^\top\otimes I_N^{\otimes 2})(\bar B^\circ)^k(e_j\otimes u)\|_2,
\]
using Lemma \ref{lm:cone} (and its straightforward extension) and noting all vectors above are in $\cC^M$. Recalling the equivalence built in Lemma \ref{lm:expectedrate}, all the terms on the r.h.s.\ can be rewritten in the form
\begin{equation}
  \label{eq:dissectforU}
  \|\mE\left(1_{\{X_k=i\}}A_k\cdots A_1 1_{\{X_0=j\}}\right)^{\otimes 2}u\|_2,
\end{equation}
the scheme appearing in the proof of Lemma \ref{lm:radius_restricted}, i.e.\ lower bounding the norm using $\|\mE (Y^{\otimes 2})u\|_2$, now with $Y$ acquired as a product of (dependent) matrices.

The Frobenius norm can be evaluated grouping by blocks, i.e.
\[
 \|(\bar B^\circ)^k\|_F^2 = \sum_{i,j\in M}\|(e_i^\top\otimes I_N^{\otimes 2})(\bar B^\circ)^k(e_j\otimes I_N^{\otimes 2})\|_F^2.
\]
The terms in the sum can be recast as
\[
 \|\mE\left(1_{\{X_k=i\}}A_k\cdots A_1 1_{\{X_0=j\}}\right)^{\otimes 2}\|_F^2,
\]
with the same random matrix appearing as in \eqref{eq:dissectforU}. The comparison argument of Lemma \ref{lm:radius_restricted} can be used, with now a maximum contrasted against a sum leading to an extra $M^2$ factor loss. But when taking $k\to \infty$ to recover $\rho(\bar B^\circ)$ and $\rho_U(\bar B^\circ)$, this vanishes.

Similarly as before, $\cC^M=\cC\times\cdots\times\cC$ is an invariant proper pointed cone for $\bar B^\circ,{\bar B}^{\circ\cG}$, and $\bar B^\circ-{\bar B}^{\circ\cG}$ inside $U^M$. Thus there exists $v\in\cC^M$ as the (right) principal eigenvector for ${\bar B}^{\circ\cG}$, for which we can compare the actions as
\begin{equation}
  \label{eq:Bactions}
  (\bar B^\circ)^k v = ({\bar B}^{\circ\cG} + (\bar B^\circ - {\bar B}^{\circ\cG}))^k v = ({\bar B}^{\circ\cG})^k v + z_k,
\end{equation}
with $z_k\in\cC^M$ thanks to the invariance properties. For any $a=(a_1,\ldots,a_M) \in\cC^M$ the norms relate easily as $\|a\|_2^2=\|a_1\|_2^2+\ldots+\|a_M\|_2^2$, thus applying Lemma \ref{lm:cone} coordinatewise we immediately get $\|a+b\|_2\geq \|a\|_2$ for $a,b\in\cC^M$.

In the end, the spectral radius of $\bar B^\circ$ is bounded below by the rate of $(\bar B^\circ)^k v$, itself bounded below by the spectral radius of ${\bar B}^{\circ\cG}$ thanks to the observations above, completing the proof.
\end{proof}

\section{Numerical experiments}
\label{sec:numerical}

\subsection{Comparing compute time of various methods on the same problem}

For comparison of the computational complexity of our newly defined method presented as {\bf Algorithm 1} and those previously described in \cite{gerencsr2019tight}, \cite{gkpushsumsynchr2024} we ran the same benchmark process demonstrated in \cite{gkpushsumsynchr2024} with transitive graphs of varying sizes -- in particular random Cayley graphs for $A_5,S_5$ and a grid -- and appended the results ($T(\iota)$) to {\it Table 1} from the aforementioned work. The compute times are collected in Table \ref{tab:benchmark}. 

\begin{table}[H]
	\centering
	\begin{tabular}{|r |r |r  |r |r |r|}
		\hline
		$N$ & $T(\gamma)$ [s] & $T(\eta_2)$ [s] & $T(\iota)$ [s] & $\log_N  ({T(\eta_2)}/{T(\gamma)})$ & $\log_N  ({T(\eta_2)}/T(\iota))$ \\
		\hline
		60  & 0.00015 & 0.18516  & 0.39510   &  1.733 &  -0.178 \\
		120 & 0.00033 & 9.12899    & 0.86568   &  2.13 &   0.492 \\
		210 & 0.00110 & 136.40832 & 2.08397   & 2.19 & 0.782 \\
		360 & 0.00514 & - & 10.6729 & - & -\\
		720 & 0.01845 & - &89.3844 & - & - \\
		\hline
	\end{tabular}
	\caption{Comparing runtimes[s] of the computation of $\gamma$, $\eta_2$ and $\iota$. $N$ refers to the size of the graph, $\gamma$ is our polynomial method presented in \cite{gkpushsumsynchr2024}, $\eta_2$ corresponds to the tensor-based computation provided by \cite{gerencser2022computable} and $\iota$ denotes the method defined by {\bf Algorithm 1}.}
	\label{tab:benchmark}
\end{table}

While the method corresponding to $\gamma$ \cite{gkpushsumsynchr2024} is still superior regarding computation time thanks to its proportional dimension reduction technique, which relies heavily on the strict symmetry of the underlying graph, we can see that Algorithm 1 is still genuinely faster than direct computation of the spectral radius for larger graph sizes.

\subsection{Gossip along a random walk with limited backtracking}

\begin{figure}[H]
\begin{center}
	 
	 \begin{subfigure}[t]{0.47\textwidth}
	 	
	 \includegraphics[width=\textwidth]{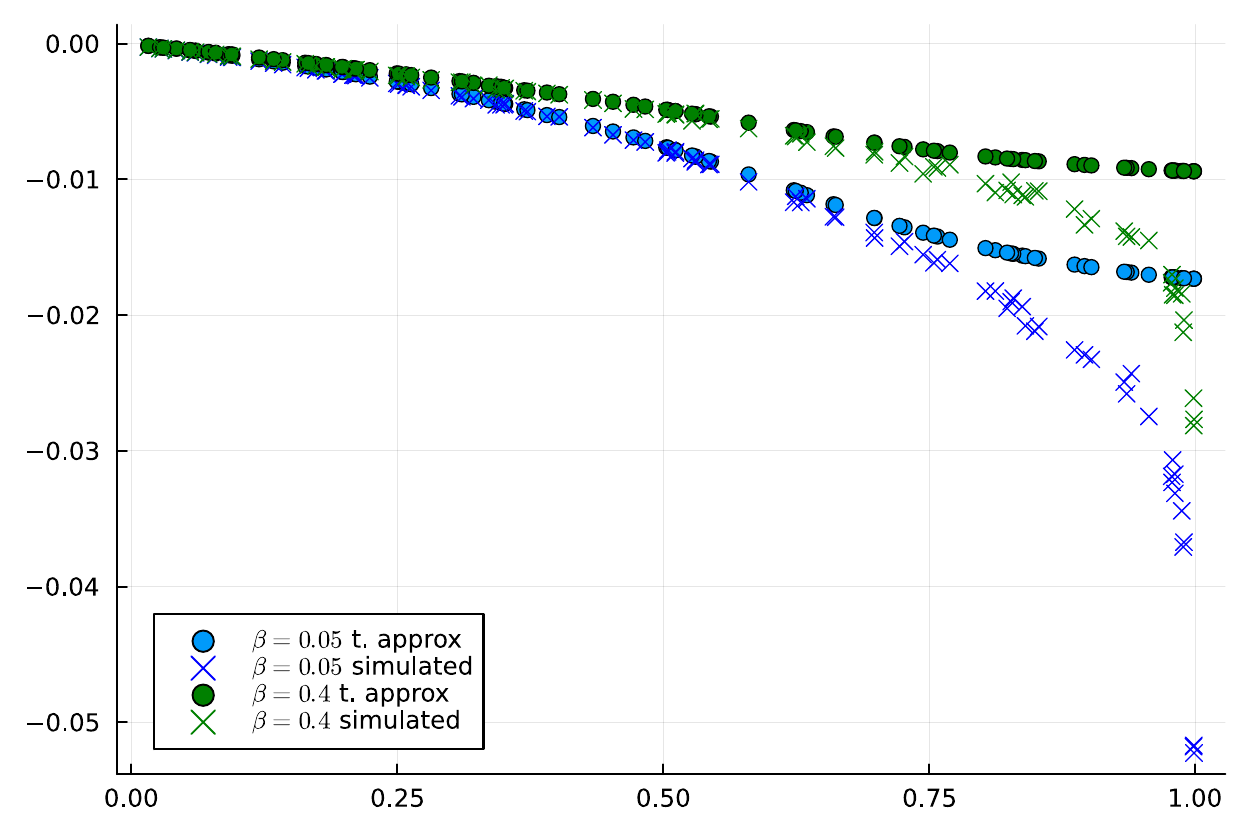}
		\caption{$5\times 6$ periodic grid}
		
\end{subfigure}
\begin{subfigure}[t]{0.47\textwidth}
	\includegraphics[width=\textwidth]{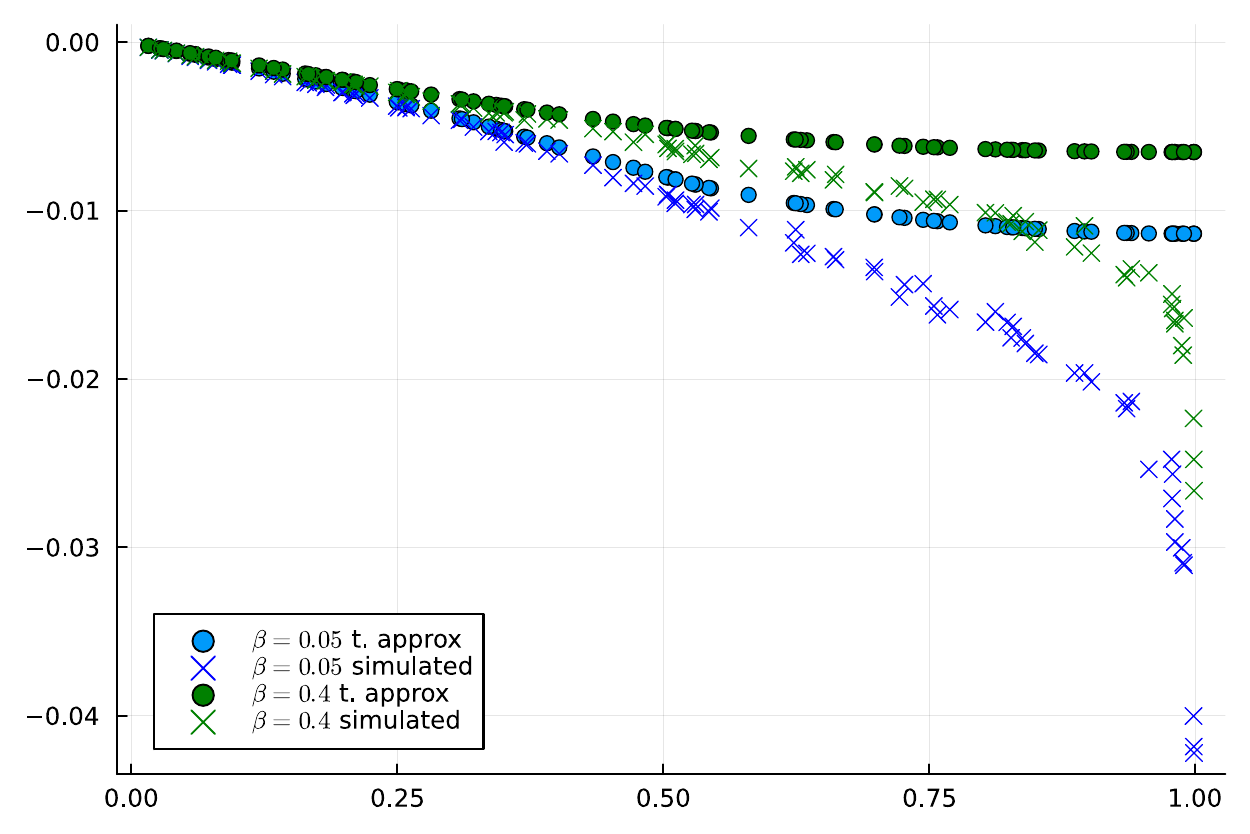}
	\caption{Random Barabási-Albert graph with $N=25$}
\end{subfigure}	 
\end{center}
\caption{Limited backtracking numerical results}
\label{fig:backtrack}
\end{figure}

In order to demonstrate the performance of our algorithm we have conducted simulations on various digraphs. A Markovian protocol was used for messaging, with state space consisting of the directed edges and transition probabilities given by 
$$ \mathbb P\big( (i,j)\ |\ (k, l)\big) = \begin{cases}
	 (1-\beta)/ (N_i-1) & \mbox{if } l=i \mbox{ and } k\neq j, \\
	 \beta  & \mbox{if } (j,i)=(k,l), \\ 
	 0 & \mbox{otherwise},
\end{cases} $$ 
where $N_i$ denotes the number of neighbors of $i$ and $\beta$ denotes the probability of backtracking. The probability of emission, i.e. $i$ sending a message to a neighbor node $j$, given the Markov state $e=(k,l)$, equals 
$$ \mathbb P\big(i\to j\ |\ (k,l)\big) = \delta_{(i,j),(k,l)}. $$
Note that this corresponds exactly to a messaging protocol driven by a random walk with limited backtracking on $\overrightarrow{G}$. Figure \ref{fig:backtrack} compares simulated convergence rates to bounds obtained by {\bf Algorithm 2a} on two separate graphs. 

Here $q$ denotes the message ratio sent by each sender.
We observed small gaps up until $q=0.5$, where these begin to widen and eventually drift apart as $q\to1$.

\subsection{Derandomization}

To illustrate Section~\ref{sec:derandomize} we used a connected Watts--Strogatz graph with parameters $N=50,\ k=12,\ \beta=0.8 $, where at each step a uniformly random vertex $v_t$ pushes to a uniformly random neighbor $u_t$. 
\begin{figure}[H]
	\centering
	\begin{subfigure}[t]{0.47\textwidth}
		\includegraphics[width=0.9\linewidth]{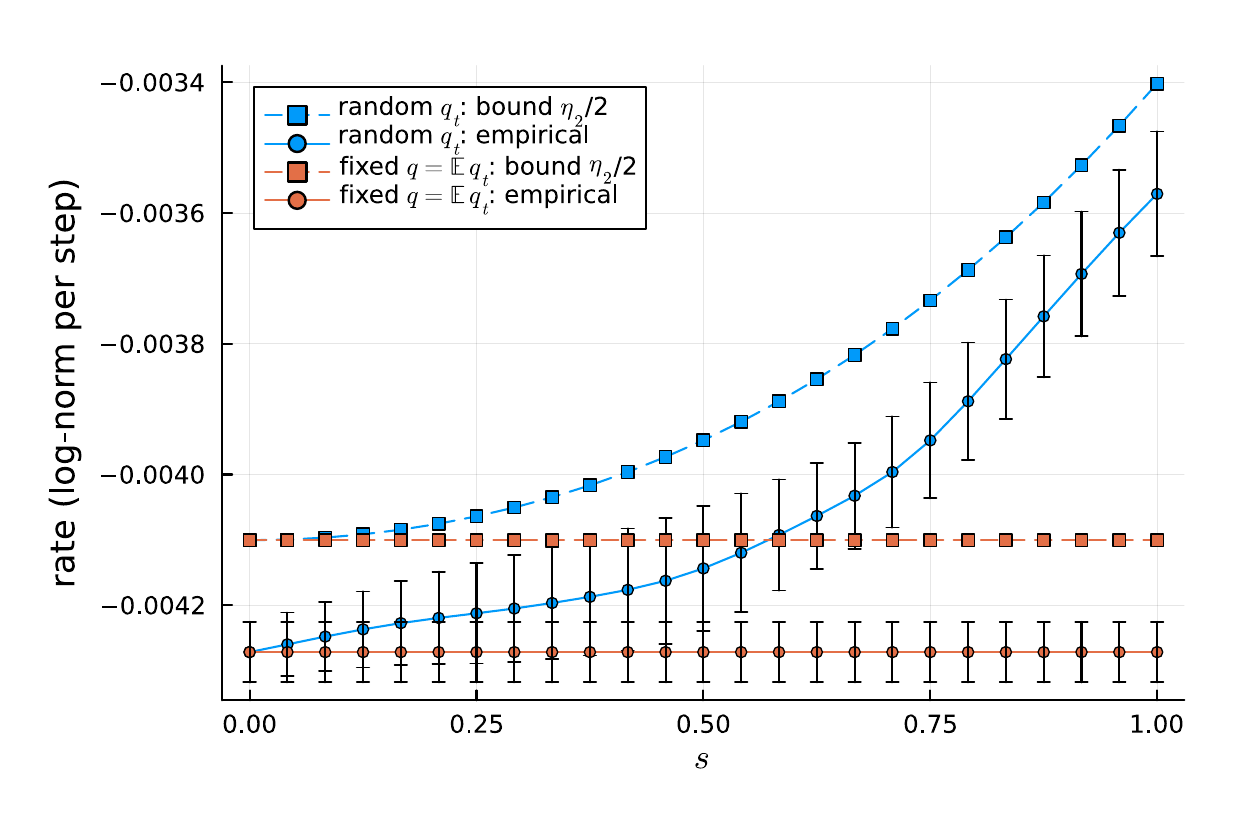}
		\caption{Example~1, as a function of the spread $s$.}
		\label{fig:derandom_ex1}
	\end{subfigure}
	\begin{subfigure}[t]{0.47\textwidth}
		\includegraphics[width=0.9\linewidth]{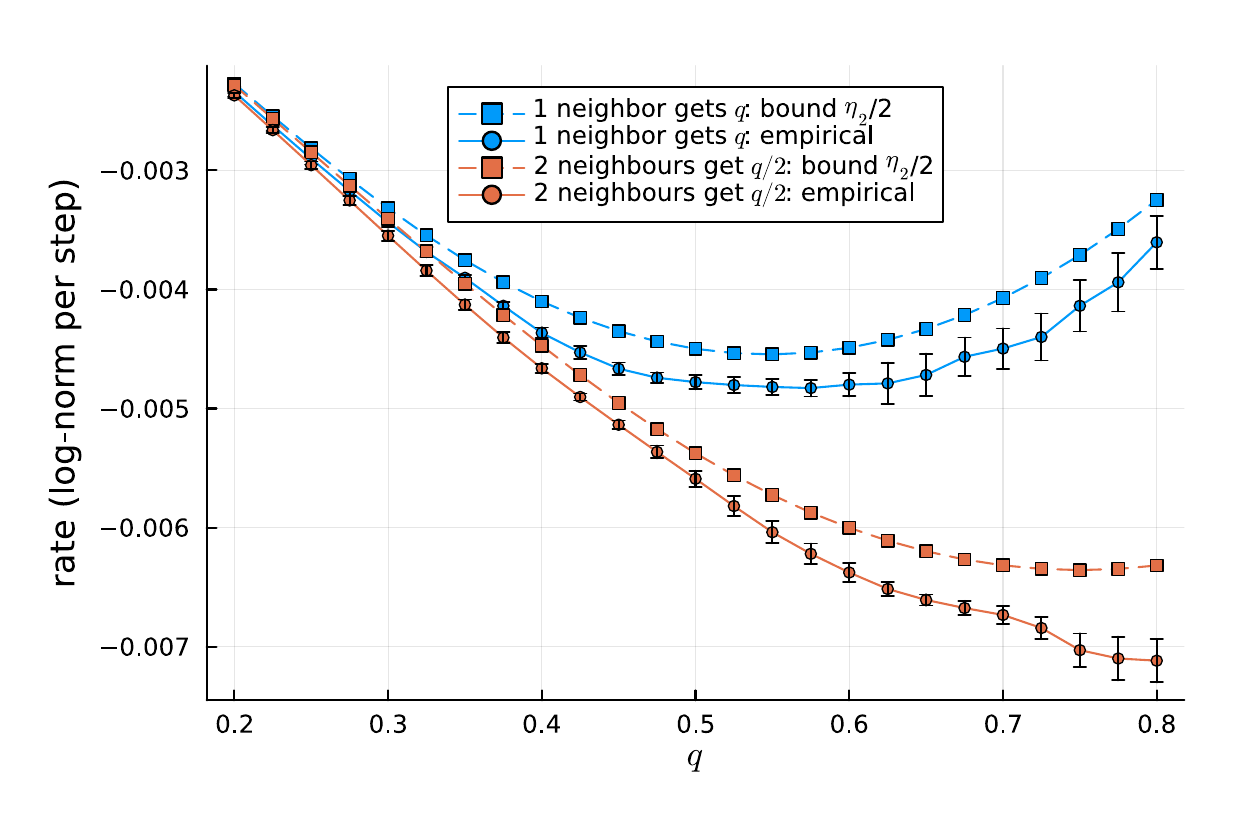}
		\caption{Example~2, as a function of $q$.}
		\label{fig:derandom_ex2}
	\end{subfigure} 
	\caption{Rates and bounds for derandomization examples}
	\label{fig:derandom}
\end{figure}
 In Example~1 the proportion $q_t$ is uniform on $[q(1-s),q(1+s)]$ with $q=0.4$, and is compared with the fixed proportion $\mE q_t=q$. The parameter $s$ controls the spread of the distribution.\\ In Example~2 the proportion $q$ is fixed, and sending it to one of two random neighbors chosen by a coin is compared with sending $q/2$ to both. We computed $\eta_2/2$ exactly from $\mE A_1^{\otimes 2}$, and also the empirical rate of $\|\mathrm{diag}(1/w_t)A_t\ldots A_1(I-J)\|_F$ (slope of its logarithm between $t/2$ and $t$, $t=20\,000$, multiple precision, mean and standard error of $10$ runs sharing the sequence $(u_t,v_t)$ within each pair). The numerical results are presented in Figure~\ref{fig:derandom}, where $\eta_2$ is the upper bound obtained from the tensor estimation method.

Example~1 on panel (a) of Figure \ref{fig:derandom} shows that randomizing $q$ slows the convergence rate, which can be explained by the decreasing trend of typical out-flux per vertex on the long run.

Example~2 on panel (b) is numerical evidence for the fact that sending to two neighbors half the prescribed quantity is asymptotically more beneficial compared to the alternative: sending the full package to one neighbor. The same pattern emerges, although in a slightly different setup, in \cite{gerencser2022computable}.

\section{Conclusions}
\label{sec:conclusions}
We have presented tools to bound the almost sure exponential convergence rate of push-sum consensus algorithms driven by HMM message passing, also tools for efficient computation of this bound, and finally tools to propose variants of the consensus algorithm for improved rates.

Simulation results demonstrate the practical usability of the analytic findings. In particular, we have observed how the convergence rate of message passing with less backtracking is better and it is approximated closely by our bounding quantity. Note that this process is conceptually close to One-Way Averaging \cite{benezit2010weighted}, for which fast consensus has been observed but there are no strong bounds on the rate yet.

Multiple open questions arise as natural follow-up directions. We have demonstrated a transformation of the protocol to improve the bound on the rate, a more flexible class of transformations would be beneficial that is rich enough to perform optimization on the protocol. We have put strong emphasis on computational accessibility, later work could identify more subtle bottlenecks besides raw size. Going one step further, our results for convergence rates are inherently in the asymptotic regime; finding bounds for the initial transient regime is a sought-after complementary study. 
\section{Declaration of Generative AI and AI-Assisted Technologies}
The core ideas and concepts as well as the proofs, serving as the backbone of this work, are intellectual products of the authors solely. The first version of the whole manuscript is human written, AI (Claude Sonnet 5.5 and Opus 5.5) was used for minor fine-tuning, improving and enhancement of writing style and prose, as well as for grammar checks. Additionally, proofs were also checked mathematically and commented on by these models, and in case any recommendations were made, the authors inspected them and applied corrections and improvements manually whenever they were well based. \\
The same models were used to assist writing the code for the simulations presented in Section \ref{sec:numerical}. Results were manually checked by the authors.

\section*{Acknowledgements}
We thank Rapha\"el Jungers for his inspiring remarks. 

\bibliographystyle{siam}
\bibliography{generic,pushsum}

\begin{thebibliography}{10}

\bibitem{benezit2010weighted}
{\sc F.~B{\'e}n{\'e}zit, V.~Blondel, P.~Thiran, J.~N. Tsitsiklis, and
  M.~Vetterli}, {\em Weighted gossip: Distributed averaging using non-doubly
  stochastic matrices}, in {P}roceedings of 2010 {IEEE} {I}nternational
  {S}ymposium on {I}nformation {T}heory ({ISIT}), 2010, pp.~1753--1757.

\bibitem{diaconis1999iterated}
{\sc P.~Diaconis and D.~Freedman}, {\em Iterated random functions}, SIAM
  review, 41 (1999), pp.~45--76.

\bibitem{gerencser2022computable}
{\sc B.~Gerencs{\'e}r}, {\em Computable convergence rate bound for ratio
  consensus algorithms}, IEEE Control Systems Letters, 6 (2022),
  pp.~3307--3312.

\bibitem{gerencsr2019tight}
{\sc B.~Gerencs{\'e}r and L.~Gerencs{\'e}r}, {\em Tight bounds on the
  convergence rate of generalized ratio consensus algorithms}, IEEE
  Transactions on Automatic Control, 67 (2022), pp.~1669--1684.

\bibitem{gerencser2018push}
{\sc B.~Gerencs{\'e}r and J.~M. Hendrickx}, {\em Push sum with transmission
  failures}, IEEE Transactions on Automatic Control, 64 (2018), pp.~1019--1033.

\bibitem{gkpushsumpolygen2025}
{\sc B.~Gerencs{\'e}r and M.~Kornyik}, {\em Low complexity convergence rate
  bounds for push-sum algorithms with homogeneous correlation structure}, arXiv
  preprint arXiv:2507.16601,  (2025).

\bibitem{gerencser2008lq}
{\sc L.~Gerencs{\'e}r and Z.~Orlovits}, {\em Lq-stability of products of
  block-triangular stationary random matrices}, Acta Scientiarum
  Mathematicarum, 74 (2008), pp.~927--944.

\bibitem{gkpushsumsynchr2024}
{\sc B.~Gerencsér and M.~Kornyik}, {\em Low complexity convergence rate bounds
  for the synchronous gossip subclass of push-sum algorithms}, IEEE Control
  Systems Letters, 8 (2024), pp.~1283--1288.

\bibitem{haagerup1985grothendieck}
{\sc U.~Haagerup}, {\em The {G}rothendieck inequality for bilinear forms on
  {C*}-algebras}, Advances in Mathematics, 56 (1985), pp.~93--116.

\bibitem{hadjicostis2014average}
{\sc C.~N. Hadjicostis and T.~Charalambous}, {\em Average consensus in the
  presence of delays in directed graph topologies}, IEEE Transactions on
  Automatic Control, 59 (2014), pp.~763--768.

\bibitem{hadjicostis2018distributed}
{\sc C.~N. Hadjicostis, A.~D. Dom{\'\i}nguez-Garc{\'\i}a, and T.~Charalambous},
  {\em Distributed averaging and balancing in network systems}, Foundations and
  Trends in Networking, 5 (2018), pp.~99--292.

\bibitem{hadjicostis2015robust}
{\sc C.~N. Hadjicostis, N.~H. Vaidya, and A.~D. Dom{\'\i}nguez-Garc{\'\i}a},
  {\em Robust distributed average consensus via exchange of running sums}, IEEE
  Transactions on Automatic Control, 61 (2015), pp.~1492--1507.

\bibitem{iutzeler2013analysis}
{\sc F.~Iutzeler, P.~Ciblat, and W.~Hachem}, {\em Analysis of sum-weight-like
  algorithms for averaging in wireless sensor networks}, IEEE Transactions on
  Signal Processing, 61 (2013), pp.~2802--2814.

\bibitem{kempe2003gossip}
{\sc D.~Kempe, A.~Dobra, and J.~Gehrke}, {\em Gossip-based computation of
  aggregate information}, in Proceedings of 44th Annual IEEE Symposium on
  Foundations of Computer Science, 2003, pp.~482--491.

\bibitem{lin2022subgradient}
{\sc Y.~Lin and J.~Liu}, {\em Subgradient-push is of the optimal convergence
  rate}, in 2022 IEEE 61st Conference on Decision and Control (CDC), IEEE,
  2022, pp.~5849--5856.

\bibitem{nedic2014distributed}
{\sc A.~Nedi{\'c} and A.~Olshevsky}, {\em Distributed optimization over
  time-varying directed graphs}, IEEE Transactions on Automatic Control, 60
  (2014), pp.~601--615.

\bibitem{nedic2017achieving}
{\sc A.~Nedi{\'c}, A.~Olshevsky, and W.~Shi}, {\em Achieving geometric
  convergence for distributed optimization over time-varying graphs}, SIAM
  Journal on Optimization, 27 (2017), pp.~2597--2633.

\bibitem{pisier2003introduction}
{\sc G.~Pisier}, {\em Introduction to operator space theory}, vol.~294,
  Cambridge University Press, 2003.

\bibitem{bahman_ps_optim2025}
{\sc P.~Rezaeinia, B.~Gharesifard, and T.~Linder}, {\em Distributed
  optimization with uncertain communications}, IEEE Transactions on Automatic
  Control, 70 (2025), pp.~2746--2753.

\bibitem{hespanha2018broadcast}
{\sc D.~Silvestre, J.~P. Hespanha, and C.~Silvestre}, {\em Broadcast and gossip
  stochastic average consensus algorithms in directed topologies}, IEEE
  Transactions on Control of Network Systems, 6 (2018), pp.~474--486.

\bibitem{tsitsiklis1997lyapunov}
{\sc J.~N. Tsitsiklis and V.~Blondel}, {\em The {L}yapunov exponent and joint
  spectral radius of pairs of matrices are hard—when not impossible—to
  compute and to approximate}, Mathematics of Control, Signals and Systems, 10
  (1997), pp.~31--40.

\end{thebibliography}

\appendix
\section{Comparison of spectral radii}
\label{sec:tensorcs}
\label{app:tensorcs}
We now show a Cauchy--Schwarz type inequality for random tensor products, analogous to the one in \cite{gerencser2008lq}. Observe that the class to be considered below is rich enough: any matrix can be represented as $\mE T_{XY}$ for appropriate (random) $T_{XY}$. Although the result can be obtained through abstract operator theoretical reasoning building on \cite{pisier2003introduction},\cite{haagerup1985grothendieck}, we think it is instructive to present an elementary and self-contained proof.

\begin{prop}
  \label{prop:tensorcs}
  Let $X_{i,j}, 1\leq i,j\leq M$ be random matrices of size $N_X\times N_X$, similarly $Y_{i,j}, 1\leq i,j\leq M$ of size $N_Y\times N_Y$ (with no independence imposed).
  Define the (random) block matrix $T_{XY}$ of $M\times M$ blocks, with $X_{i,j}\otimes Y_{i,j}$ in cell $(i,j)$. Analogously define $T_{XX}, T_{YY}$. Then there exists a constant $c>0$ depending only on the dimensions so that
  \[
    \|\mE T_{XY}\|^2 \leq  c   \|\mE T_{XX}\| \cdot \|\mE T_{YY}\|.
  \]
  Moreover, for the spectral radii we have
  \[
    \rho( \mE T_{XY})^2 \leq  \rho(\mE T_{XX}) \cdot \rho(\mE T_{YY}).
  \]
\end{prop}

\begin{proof}
  We will show the first inequality for a specialized norm, then by the equivalence of norms it follows for other norms by adjusting $c$, leading to the first claim.
  Let us define test vectors as
  \begin{align*}
    S_{XY}&=\{u_X\otimes u_Y ~:~ u_X\in \mR^{N_X}, u_Y\in \mR^{N_Y}, \|u_X\|=\|u_Y\| = 1\}\\
    S_{XY}^{(M)}&=\{u= (u_1^\top,\ldots,u_M^\top)^\top ~:~u_i\in S_{XY}\}.
  \end{align*}
  Analogously we can define $S_{XX},S_{XX}^{(M)},S_{YY},S_{YY}^{(M)}$, and the following reasoning also carries through for these variants. Note that $S_{XY}$ generates the full vector space $\mR^{N_X N_Y}$ and $S_{XY}^{(M)}$ generates $\mR^{MN_XN_Y}$.
  Now for a matrix of size $\mR^{MN_XN_Y\times MN_XN_Y}$ we can set our norm as follows:
  \begin{equation}
    \label{eq:starnorm}
    \|T\|_* = \sup\left\{u^\top T v \mid u,v \in S_{XY}^{(M)}\right\}.
  \end{equation}
  Note that the supremum is actually a maximum as we optimize over a compact space.
  We need to confirm this definition is indeed a norm, to be able to rely on the equivalence of norms afterwards. 

  It scales linearly directly from the definition. In order to show that only the 0 matrix has norm 0, consider any $0\neq T\in\mR^{MN_XN_Y\times MN_XN_Y}$. Then there exists $\tilde u,\tilde v\in\mR^{MN_XN_Y}$ so that $\tilde u^\top T \tilde v>0$. As $\tilde u,\tilde v$ are both linear combinations of vectors in $S_{XY}^{(M)}$, there must be $u,v\in S_{XY}^{(M)}$ so that $u^\top T v \neq 0$. This confirms $\|T\|_*>0$ (testing with $u,v$ or $u,-v$). For the triangle inequality, for any matrices $T,T'\in\mR^{MN_XN_Y\times MN_XN_Y}$ and vectors $u,v\in S_{XY}^{(M)}$ the linearity
   of the quadratic forms $u^\top (T+T') v=u^\top T v+u^\top T' v$ holds.

Taking the supremum in $u,v$ for the l.h.s.\ expresses $\|T+T'\|_*$, these vectors might be suboptimal for the r.h.s.\ terms, separately further optimizing them leads to $\|T+T'\|_*\leq\|T\|_*+\|T'\|_*$.

 Turning to the actual claim, for the case of $T_{XY}$ defined in the statement, a block matrix of tensor products, the quadratic form appearing in the definition \eqref{eq:starnorm} of $\|\cdot\|_*$ can be rewritten as
 \begin{align}
   \label{eq:norm_quadratic}
   u^\top T_{XY} v &= \sum_{i,j=1}^M \left(u_{X,i}^\top X_{i,j} v_{X,j}\right)\cdot \left(u_{Y,i}^\top Y_{i,j} v_{Y,j}\right)
 \end{align}
 Taking expectations and using the Cauchy--Schwarz inequality both for the sum of products and then for the expectation of products we get
 \begin{align*}
   u^\top \mE T_{XY} v = \mE (u^\top T_{XY} v)   &\leq \mE\left(\sqrt{\sum_{i,j=1}^M \left(u_{X,i}^\top X_{i,j} v_{X,j}\right)^2}\cdot \sqrt{\sum_{i,j=1}^M\left(u_{Y,i}^\top Y_{i,j} v_{Y,j}\right)^2}\right)\\
  &\leq \sqrt{\mE\sum_{i,j=1}^M \left(u_{X,i}^\top X_{i,j} v_{X,j}\right)^2}\cdot \sqrt{\mE\sum_{i,j=1}^M\left(u_{Y,i}^\top Y_{i,j} v_{Y,j}\right)^2}.
 \end{align*}
 Similarly as in \eqref{eq:norm_quadratic}, under the square roots we see quadratic forms for $T_{XX}$ and $T_{YY}$. Maximizing the l.h.s.\ in $u,v$ we obtain $\|\mE T_{XY}\|_*$, while on the r.h.s.\ we get lower bounds for $\sqrt{\|\mE T_{XX}\|_*}$ and $\sqrt{\|\mE T_{YY}\|_*}$, thus
 \begin{equation}
   \label{eq:tensor-cs}
   \|\mE T_{XY}\|_*^2 \leq  \|\mE T_{XX}\|_* \cdot \|\mE T_{YY}\|_*,  
 \end{equation}
  confirming the first claim relying on the equivalence of norms.

  The comparison of spectral radii will follow from the above, we will find a way to compare the norms of $(\mE T_{XY})^\ell,(\mE T_{XX})^\ell,(\mE T_{YY})^\ell$.

Assume another collection of random matrices $X'_{i,j}, Y'_{i,j}$ is given with the same (joint) distribution as $X_{i,j},Y_{i,j}$, generating $T_{X' Y'}, T_{X' X'}, T_{Y' Y'}$ further assuming they are independent of $X_{i,j},Y_{i,j}$. For the $(i,j)$ block of $T_{X' Y'}T_{X Y}$ we get
  \[
    \sum_{k=1}^M (X'_{i,k}\otimes Y'_{i,k})(X_{k,j}\otimes Y_{k,j}) = \sum_{k=1}^M (X'_{i,k}X_{k,j})\otimes (Y'_{i,k}Y_{k,j}).
  \]
  Motivated by this equation, let us define $\tilde X_{i,j},\tilde Y_{i,j}$ as follows. Let $K$ be uniformly random on $[M]$, independent of the other variables, and set
  \begin{align*}
    \tilde X_{i,j} &:= \sum_{k=1}^M \delta_k(K)X'_{i,k}X_{k,j},&
    \tilde Y_{i,j} &:= \sum_{k=1}^M \delta_k(K)Y'_{i,k}Y_{k,j}.
  \end{align*}
  For these random variables we have
  \[
    \mE T_{\tilde X \tilde Y} = \frac 1M \mE (T_{X' Y'}T_{X Y}) = \frac 1M \mE T_{X' Y'} \mE T_{X Y} = \frac 1M (\mE T_{X Y})^2.
  \]
  Iterating this construction we can multiply several block matrices, in particular we can express $(\mE T_{X Y})^\ell$, using some random matrices $X^\ell_{i,j},Y^\ell_{i,j}$ so that
  \begin{align*}
    \mE T_{X^\ell Y^\ell} &= M^{-\ell+1} (\mE T_{X Y})^\ell,\\
    \mE T_{X^\ell X^\ell} &= M^{-\ell+1} (\mE T_{X X})^\ell,\\
    \mE T_{Y^\ell Y^\ell} &= M^{-\ell+1} (\mE T_{Y Y})^\ell.
  \end{align*}

  In summary, we can write
  \begin{align*}
    \rho( \mE T_{XY})^2 = \lim_{\ell\to\infty} \|(\mE T_{X Y})^{2 \ell}\|_*^{1/\ell} &= \lim_{\ell\to\infty} (M^{2\ell-1}\|\mE T_{X^{2\ell} Y^{2\ell}}\|_*)^{1/\ell}\\ &\leq \lim_{\ell\to\infty} (M^{4\ell-2}\|\mE T_{X^{2\ell} X^{2\ell}}\|_* \cdot \|\mE T_{Y^{2\ell} Y^{2\ell}}\|_*)^{1/2\ell} =  \rho(\mE T_{XX}) \cdot \rho(\mE T_{YY}),
  \end{align*}
essentially using \eqref{eq:tensor-cs} in the middle. This confirms our second claim.

\end{proof}

The comparison needed for Algorithm 2 becomes a clean corollary.

\begin{proof}[Proof of Lemma \ref{lm:BC_rhos}]
  Applying Proposition \ref{prop:tensorcs} with $X_{i,j} = \sqrt{p_{ji}}$ and $Y_{i,j} = \sqrt{p_{ji}}A_{1,i}(I-J)$ implies $\mE T_{XY} = \bar C^{\circ}$, $\mE T_{XX} = P^\top$ and $\mE T_{YY} = \bar B^{\circ}$. Furthermore
  by noting $\rho(P)=1$ we conclude the proof.    
\end{proof}

\end{document}